\documentclass[12pt]{article}

\usepackage[utf8]{inputenc}
\usepackage[T1]{fontenc}
\usepackage{amsmath}
\usepackage{amsfonts}
\usepackage{latexsym}
\usepackage[]{graphicx}  
\usepackage{amssymb}
\usepackage{amsthm}
\usepackage{thmtools}
\usepackage[margin=2cm]{geometry}
\usepackage{url}
\usepackage{color}
\usepackage{enumerate}
\usepackage[shortlabels]{enumitem} 
\usepackage[noadjust]{cite}
\usepackage{framed}
\usepackage{ifdraft}
\usepackage{relsize}
\usepackage{caption}
\usepackage{subcaption}
\usepackage{mathtools}
\usepackage{xspace}
\usepackage[affil-it]{authblk}
\usepackage{xcolor}
\usepackage{comment}
\ifdraft{
	\usepackage[colorinlistoftodos]{todonotes}
	\usepackage{letltxmacro}
	\LetLtxMacro{\oldtodo}{\todo}
	\usepackage{xspace}
	\renewcommand{\todo}[2][]{\oldtodo[#1]{#2}\xspace}%
	\usepackage{xcolor}
	\usepackage[mark,dirty={Modified!}]{gitinfo2}
	
}{
	\usepackage[disable]{todonotes}
}

\usepackage{pdftexcmds}
\makeatletter
\newcommand{\stringcases}[3]{%
	\romannumeral
	\str@case{#1}#2{#1}{#3}\q@stop
}
\newcommand{\str@case}[3]{%
	\ifnum\pdf@strcmp{\unexpanded{#1}}{\unexpanded{#2}}=\z@
	\expandafter\@firstoftwo
	\else
	\expandafter\@secondoftwo
	\fi
	{\str@case@end{#3}}
	{\str@case{#1}}%
}
\newcommand{\str@case@end}{}%

\long\def\str@case@end#1#2\q@stop{\z@#1}
\makeatother

\usepackage[draft=false,breaklinks,backref=page,colorlinks]{hyperref} 

\renewcommand*{\backref}[1]{}
\renewcommand*{\backrefalt}[4]{%
	\ifcase #1 %
	\relax %
	\or
	(page #4).%
	\else
	(pages #4).%
	\fi%
}

\usepackage[nameinlink]{cleveref}

\newtheorem{theorem}{Theorem}[section]
\newtheorem{lemma}[theorem]{Lemma}

\newtheorem*{proposition*}{Proposition}
\newtheorem{corollary}[theorem]{Corollary}
\newtheorem*{corollary*}{Corollary}

\newtheorem*{conjecture*}{Conjecture}
\newtheorem{proposition}[theorem]{Proposition}

\theoremstyle{definition}
\newtheorem{definition}[theorem]{Definition}
\newtheorem{example}[theorem]{Example}

\newtheorem{alg}[theorem]{Algorithm}

\theoremstyle{remark}
\newtheorem{remark}[theorem]{Remark}

\numberwithin{equation}{section}

\newcommand{\N}{\mathbb{N}}
\newcommand{\Z}{\mathbb{Z}}
\newcommand{\Q}{\mathbb{Q}}
\newcommand{\R}{\mathbb{R}}
\newcommand{\C}{\mathbb{C}}

\newcommand{\Id}{{\mathrm{Id}}}

\newcommand{\A}{\mathcal A}

\newcommand{\Abb}{\mathbb{A}}
\newcommand{\Ocal}{\mathcal{O}}

\newcommand{\SL}{\mathrm{SL}}

\newcommand{\mat}[4]{\begin{pmatrix}
		#1 & #2 \\ #3 & #4
\end{pmatrix}}

\newcommand{\vek}[2]{\begin{pmatrix}#1\\#2\end{pmatrix}}
\newcommand{\tvek}[3]{\begin{pmatrix}#1\\#2\\#3\end{pmatrix}}

\crefname{theorem}{theorem}{theorems}
\crefname{corollary}{corollary}{corollaries}
\crefname{example}{example}{examples}
\crefname{lemma}{lemma}{lemmas}
\crefname{proposition}{proposition}{propositions}
\crefname{definition}{definition}{definitions}
\crefname{observation}{observation}{observations}

\begin{document}

	\providecommand{\keywords}[1]
	{
		\small	
		\textbf{\textit{Keywords:}} #1
	}

		\title{Periodicity of general multidimensional continued fractions using repetend matrix form}

	\author{Hanka \v Rada%
		\thanks{Present adress: Czech Technical University in Prague, Faculty of Information Technology, Department of Applied Mathematics, Th\' akurova 9, 160 00 Prague 6, Czech republic
			\\Electronic address: \texttt{hanka.rada@fit.cvut.cz}}}
	\affil{\footnotesize Czech Technical University in Prague\\ Faculty of Nuclear Sciences and Physical Engineering\\ B\v rehová 78/7, 115 19 Star\' e M\v esto\\ Czech Republic}
	
	\author{\v St\v ep\' an Starosta%
		\thanks{Electronic address: \texttt{stepan.starosta@fit.cvut.cz}}}
	\affil{\footnotesize Czech Technical University in Prague\\	Faculty of Information Technology\\ Department of Applied Mathematics\\ Th\' akurova 9, 160 00 Prague 6\\ Czech republic}

\author{V\' \i t\v ezslav Kala%
	\thanks{Electronic address: \texttt{vitezslav.kala@matfyz.cuni.cz}}}
\affil{\footnotesize Charles University\\ Faculty of Mathematics and Physics\\ Department of Algebra\\ Sokolov\-sk\' a 83, 18600 Praha~8\\ Czech Republic}

	\date{}
	
	\maketitle

\begin{abstract}
We consider expansions of vectors by a general class of multidimensional continued fraction algorithms. If the expansion is eventually periodic, then we describe the possible structure of a matrix corresponding to the repetend and use it to prove that a number of vectors have an eventually periodic expansion in the Algebraic Jacobi--Perron algorithm. Further, we give criteria for vectors to have purely periodic expansions; in particular, the vector cannot be totally positive.
\end{abstract}
\maketitle

\section{Introduction}

Lagrange showed that any quadratic irrational number has an eventually periodic classic continued fraction, thus providing a characterization for quadratic irrational numbers in terms of their continued fraction expansions.
In 1839, Hermite opened the question of similar characterization for higher order irrationalities.
More specifically, he asked for a representation by (eventually) periodic sequences that would capture the algebraicity of the represented number, especially focusing on cubic numbers.
His question was published later in 1850, see \cite[p. 286]{Hermite1850_2}.

A multitude of algorithms generalizing the classic continued fraction algorithm exist, in general, they are called multidimensional continued fraction (MCF) algorithms.
They have various properties related to the algebraicity of the represented number and their approximation capabilities.
However, till today, there is no satisfactory answer to Hermite's question in full generality.
Allowing some constraints, partial answers are known, for instance, Murru \cite{murru2015periodic} exhibits an algorithm which provides an eventually periodic representation for any cubic irrational number provided its minimal polynomial is explicitly known.

The history of MCF algorithms starts with Jacobi in 1868, \cite{jacobi}.
Jacobi's algorithm was further generalized by Perron, \cite{Perron,perron1935satz}, forming a method known today as Jacobi--Perron algorithm.
Other well-known algorithms are Poincaré algorithm (1884, \cite{poincare}), Brun algorithm (1920, \cite{brun1920}), Selmer algorithm (1961, \cite{selmer}), and Fully subtractive algorithm (1995, \cite{schweiger1995}).

There exists a variety of MCF algorithms that further generalized the common concept, which is discussed for instance in \cite{brentjes,Bryuno_1997,schweiger2000}.
Let us mention some of these algorithms: the Modified Jacobi--Perron algorithm (Bernstein, 1965, \cite{bernstein1965rational}), the Heuristic APD-algorithm (Karpenkov, 2022, \cite{karpenkov2021hermite}), the $\rm{sin}^2$-algorithm (Karpenkov 2021, \cite{karpenkov2021periodic}), the algorithm of Abrate, Barbero, Cerruti and Murru (2013, \cite{abrate2013periodic}), and the Algebraic Jacobi--Perron algorithm (Tamura and Yasutomi,\cite{tamura2009new,tamura2012some}).
In~\cite{karpenkov2021periodic}, the author of the $\rm{sin}^2$-algorithm proves that his algorithm is periodic for every totally real cubic vector.
In~\cite{tamura2009new}, the authors of the Algebraic Jacobi--Perron algorithm provide some evidence that their algorithm may answer Hermite's question for cubic numbers; however, no proof is provided.

Besides the very general Hermite's questions, another recent source of interest in multidimensional continued fractions came from the study of universal quadratic forms over number fields. In degree two, there were numerous recent results bounding ranks of universal forms in terms of coefficients of the continued fraction for $\sqrt D$, e.g., \cite{blomer2017rank,kala2016universal}. While there were some extensions of such results also to higher degrees \cite{kala2020universal,kala2021lifting,yatsyna2019lower}, they remain quite limited. Among the motivations and hopes for the present article are that the general methods developed here may be useful also for the application to universal forms.

In addition to general results on various MCF algorithms, there are many results for specific classes of algebraic numbers, for example, by Bernstein \cite{bernstein1964periodical,bernstein1964periodicity,bernstein1964periodische}, Raju \cite{raju1976periodic}, Levesque \cite{levesque1976jacobi,levesque1977class}, Dubois and Paysant-Le Roux \cite{paysant1975nouvelles,dubois1975algorithme}, Greiter \cite{greiter1977mehrdimensionale}, Bouhamza \cite{bouhamza1984algorithme,bouhamza1984jacobi}.

All the above mentioned algorithms are similar to the Jacobi--Perron generalization of continued fraction algorithm; they provide a representation of a vector $\vec{v}$ in the form of a sequence of matrices from a fixed set.
A different approach to the generalization of classic continued fractions is provided by the so-called geometric MCF algorithms described for example in the book of Karpenkov \cite{karpenkov2013geometry}.

In this article, we focus solely on MCF algorithms of Jacobi--Perron type.
In the case when the representation given by such an algorithm is eventually periodic, we consider its repetend also in the form of a single matrix $M$ which equals to the product of the matrices of the repetend.
Under the assumption that $\vec{v}$ forms a basis of $\Q(\vec{v})$, we show (\Cref{thm:candidates_onMP_algebraic_new}) that $M$ is equal to the transposed matrix in the basis $\vec{v}$ of multiplication by an algebraic unit.
In~\Cref{th:weak_exists_Q}, we show that such a matrix can in fact be determined by its single column by a mapping depending on $\vec{v}$ only, i.e., independent on the algebraic unit.
We give an explicit form of this relation in~\Cref{prop:T_epsilon} based on the minimal polynomial of $y$ in the case that $\vec{v}$ is the polynomial basis $\left( y^{n-1}, \dots, y, 1 \right)^T$.
As an example, we refine the obtained formulas for dimension $3$ in~\Cref{prop:T_epsilon_cubic_case}.

We also study purely periodic expansions.
We give two necessary conditions on $\vec{v}$ to have a purely periodic expansion in a given MCF algorithm.
The first condition, \Cref{prop:purly_per}, states that the vector with a purely periodic expansion is not totally positive.
The second condition, \Cref{cor:not_purely_per}, states that the norm of $y$ needs to be $(-1)^{n-1}$ in order for $\vec{v} = \left( y^{n-1}, \dots, y, 1 \right)^T$ to have a purely periodic expansion.

Finally, in \Cref{sec:candidates}, we give a procedure how to find candidates on the product of the matrices of the repetend.
In fact, this procedure can be used to find the expansion in the given MCF algorithm in some cases, hence we refer to it as the repetend matrix form of the algorithm.
We use this form in \Cref{sec:application} to calculate the MCF expansions in order to show that a large class of vectors has an eventually periodic expansion in the case of Algebraic Jacobi--Perron algorithm in dimension $3$ (\Cref{thm:AJPA_periodic_class}).
It is notable that this result generalizes a previous result of Tamura and Yasutomi~\cite[Theorem 2.4]{tamura2009new}.

The article is organized as follows.
In \Cref{sec:prelim}, we give some necessary notations and introduce the MCF algorithms.
\Cref{sec:wekla_conv} is dedicated to matrices of multiplication in $\Q(\alpha)$.

\Cref{sec:periodic_exp} studies eventually periodic MCF expansions and the matrices that represent the product of their repetend.
\Cref{sec:candidates} elaborates a procedure on how to find a candidate on the product of repetend.
\Cref{sec:application} applies the previous results to prove that a large class of vectors has eventually periodic expansion in a given setting.

	\section{Preliminaries} \label{sec:prelim}

	A number $\alpha \in \C$ is \emph{algebraic} over $\Q$ if it is a root of some polynomial $f$ over $\Q$.
	The set of algebraic numbers (over $\Q$) is denoted by $\mathbb{A}$.
	If $\alpha$ and $\alpha'$ are roots of the same irreducible polynomial $f$, then we say that $\alpha'$ is a \emph{conjugate} of $\alpha$.

	The \emph{degree} of $\alpha$ is the least number $n$ such that $\alpha$ is a root of a polynomial over $\Q$ of degree $n$.
	Algebraic numbers of degree two are called \emph{quadratic} and algebraic numbers of degree three are called \emph{cubic} (they are roots of \emph{quadratic} respectively \emph{cubic} polynomial with rational coefficients).

	Let $\alpha_1,\dots,\alpha_n \in \Abb$. The \emph{number field $K=\Q(\alpha_1,\dots,\alpha_n)$} is defined by
	\[
	K=\Q(\alpha_1,\dots,\alpha_n) := \bigcap \{T|T \text{ is a subfield of } \mathbb{C}, \alpha_1,\dots, \alpha_n \in T\}.
	\]
	The \emph{degree} of the number field $K$ is the dimension of $K$ as a vector space over $\Q$. 
 
The well-known \textit{primitive element theorem} says
that for every $\alpha_1,\dots, \alpha_n \in \Abb$, there exists $\alpha \in \Abb$ such that $\Q(\alpha_1,\dots,\alpha_n) = \Q(\alpha)$. If $\alpha$ is an algebraic number of degree $n$, then
\[
\Q(\alpha) = \{ a_0+a_1 \alpha + \dots + a_{n-1} \alpha^{n-1} |a_i \in \Q \}.
\]

	A number $\beta \in \C$ is called an \emph{algebraic integer} if there is a monic polynomial $f$ over $\Z$ such that $f(\beta)=0$. The set of all algebraic integers is denoted by $\mathbb{B}$. \emph{The ring of integers} of the number field $\Q(\alpha)$ is the set $\Ocal_{\Q(\alpha)} := \Q(\alpha) \cap \mathbb{B}$.

	Let  $s: \Q(\alpha) \to \Q(\alpha)$ be a linear transformation. Moreover, let $S^B \in \Q^{n,n}$ be the matrix of the transformation $s$ in a basis $B$.
	If $S^{B_1}$ and $S^{B_2}$ are two matrices of the same transformations but in different bases, then $S^{B_1}$ is similar to $S^{B_2}$ (i.e., there exists an invertible matrix $U$ such that $S^{B_1} = US^{B_2}U^{-1}$). In particular, they have the same determinant. This means that we can define \emph{the determinant} of the transformation $s$ as $\det(s) = \det(S)$, where $S$ is an arbitrary matrix of the transformation $s$.
	
	We assign to each element $\delta \in \Q(\alpha)$ a linear transformation $t_\delta: \Q(\alpha) \to \Q(\alpha)$ which is defined by
	\begin{equation} \label{eq:t_alpha}
	t_\delta(x) = \delta x
	\end{equation}
	for every $x \in \Q(\alpha)$.
	The matrix of this transformation is denoted $T_\delta$.
	
	Let $\beta \in \mathbb{A}$ and $\gamma \in \Q(\beta)$. Then  the \emph{norm $N_{\Q(\beta)|\Q}(\gamma)$} (or simply $N(\gamma)$ if it is clear that $\gamma \in \Q(\beta)$) of $\gamma$ is the determinant of the linear transformation $t_\gamma$. In other words
	\[
	N_{\Q(\beta)|\Q}(\gamma) = \det(T_\gamma) \in \Q.
	\]

	A \emph{unit} in a ring $R$ with identity $1_R$ is an invertible element $u$ of $R$, i.e., there exists an element $v \in R$ such that $uv = vu = 1_R$. The units of a ring $R$ form a group with respect to multiplication, we call it the \emph{group of units $U(R)$} of $R$. In the ring of integers $\Ocal_{\Q(\alpha)}$ of a number field $\Q(\alpha)$, we can characterize the group of units in the following way.
	If $\beta \in \Ocal_{\Q(\alpha)}$, then $\beta \in U(\Ocal_{\Q(\alpha)})$ if and only if $N(\beta)=\pm 1$.
	Due to the Dirichlet's unit theorem, we can also determine the rank (the number of multiplicatively independent generators) of the group of units $U(\Ocal_{\Q(\alpha)})$.
	
	\begin{theorem}[Dirichlet's unit theorem] \label{thm:dirichlet}
		Let $K = \Q(\alpha)$ be a number field. The group of units of $\Ocal_K$ is finitely generated and its rank is equal to
		\[
		r = r_1+r_2-1,
		\]
		where $r_1$ is the number of real conjugates of $\alpha$ and $2r_2$ is the number of nonreal complex conjugates of $\alpha$.
	\end{theorem}
	For example, if $\alpha$ is a cubic number, then the group of units $U(\Ocal_K)$ has rank either $2$ or $1$.

	Let $r$ be the rank of $U(\Ocal_K)$. The set of units $u_1, \dots ,u_r$ is called \emph{the set of fundamental units} if it is multiplicatively independent and it generates (modulo roots of unity) the group $U(\Ocal_K)$, i.e. if every unit $u$ can be written uniquely in the form
	\begin{equation} \label{eq:unit}
	u = \zeta u_1^{m_1}\dots u_r^{m_r},
	\end{equation}
	where $m_i \in \Z$ for all $i \in \{1,\dots,r\}$ and $\zeta$ is some root of unity (i.e. there exists $p \in \Z_+$ such that $\zeta^p =1$).
	
		If $K = \Q(\alpha)$ is an algebraic number field of odd degree, then the roots of unity have the following simple form.

	\begin{theorem}[Theorem 13.5.2 in \cite{alaca2004introductory}] \label{thm:roots_of_unity}
		Let $K = \Q(\alpha)$ be an algebraic number field of odd degree. The roots of unity in $\Ocal_K$ are $\pm 1$.
	\end{theorem}

	The following subsections are dedicated to Jacobi--Perron type MCF algorithms and their elementary properties.
	
	\subsection{Jacobi--Perron type MCFs}

	Let $n$ be a positive integer.
	A Jacobi--Perron type MCF acts on $\R_+^{n}$ and it is specified by two sets, $\mathcal{I}$ and $\mathcal{A}$.
	The first set is an at most countable set of pairwise disjoint subsets of $\R_{+}^{n}$:
	\[
	\mathcal{I} = \{I_1,I_2,\dots \}
	\]
	where $\forall \alpha >0, \forall I \in \mathcal{I}, \alpha I \subseteq I$,
	while the second set is a set of invertible matrices from $\R^{n,n}$:
	\[
	\A = \{A_1,A_2,\dots\}
	 \] 
	 having the same cardinality as $\mathcal{I}$.
	 Moreover, we assume that $\A$ is such that no product of matrices from $\A$ equals the identity matrix.
	 Given these two sets, a representation of a vector $\vec{v} \in \R_+^{n}$ is obtained by the following algorithm.
	
	\begin{alg}[Multidimensional continued fraction algorithm with sets $(\mathcal{I},\mathcal{A})$] \label{alg:MCF_alg}
		Let $\vec{v} \in \R_+^{n}$.
		
		Set $\vec{v}^{(0)} \coloneqq \vec{v}, i\coloneqq 0$.
		
		Repeat:
		
		Let $j$ be some index such that $\vec{v}^{(i)} \in I_j$.
		If there is no such $j$, the algorithm stops.
		Otherwise set
		\[ 
		\vec{v}^{(i+1)} \coloneqq A_j^{-1}\vec{v}^{(i)}
		\]
		and $A^{(i)}\coloneqq A_j$.
		Set $i \coloneqq i+1$.	
	\end{alg}

	\begin{definition} \label{def:MCF}
		The sequence $(A^{(i)})_{i=0}^{\infty}$ from \Cref{alg:MCF_alg} is called \emph{an $(\mathcal{I},\A)$ $(n-1)$-dimensional continued fraction expansion} of the vector $\vec{v}$.
	\end{definition}

	If not ambiguous, we will often say only expansion of $\vec{v}$.
Moreover, we identify the expansion of $\vec{v}$ with $\vec{v}$, i.e., we write $\vec{v} = ( A^{(0)} , A^{(1)}, \dots )$.

	\begin{remark} \label{rem:renormalisation}
		Let $\alpha \in \R_+$. 
		Since the elements of $\mathcal{I}$ satisfy $\forall I \in \mathcal{I}, \alpha I \subseteq I$,
		we conclude that MCF expansions of $\vec{v}$ and of $\alpha \vec{v}$ are identical.
	\end{remark}

	In what follows, we will often use the last remark and we will work with vectors $\vec{z}^{(i)} = \lambda \vec{v}^{(i)}$, where $\lambda$ is such that $(\vec{z}^{(i)})_n = 1$, instead of $\vec{v}^{(i)}$.	
	In other words, the \Cref{alg:MCF_alg} works in a projective space, which is the reason for calling the algorithm $(n-1)$-dimensional and not $n$-dimensional as the value $n-1$ corresponds to the dimension of the underlying projective space.
Nevertheless, we prefer to work exclusively with homogeneous coordinates.
	
A MCF algorithm is \emph{unimodular} if the matrices from $\A$ are unimodular, that is, they have determinant equal to $\pm 1$.

An expansion of a vector $\vec{v}  = ( A^{(0)} , A^{(1)}, \dots )$ is \emph{eventually periodic} if there exists $N$ and positive $p$ 
such that $\vec{v}^{(N)} = \lambda \vec{v}^{(N+p)}$  for some real $\lambda$.
If follows that $A^{(i)} = A^{(i+p)}$ for all $i \geq N$.
	We write also
	\[
	\vec{v} = \left ( A^{(0)} , A^{(1)}, \dots, A^{(N-1)},   \overline{A^{(N)}, A^{(N+1)}, \dots, A^{(N+p-1)}}  \right).
	\]
	If $N = 0$, then the expansion is \emph{purely periodic}.
	
	The sequence of matrices $\left( A^{(0)} , A^{(1)}, \dots, A^{(N-1)} \right) $ is called a \emph{preperiodic part} and the sequence of matrices $\left( A^{(N)}, A^{(N+1)}, \dots, A^{(N+p-1)} \right) $ is called a \emph{repetend}.
	The number $N$ is called a \emph{preperiod} and the number $p$ is called a \emph{period}.

	It follows from \Cref{alg:MCF_alg} that
	\[
	A^{(0)} \cdots A^{(i-2)}A^{(i-1)} \vec{v}^{(i)} = \vec{v}^{(0)}
	\]
	and therefore, we shall consider the preperiodic part of an expansion and its repetend as matrices, i.e., $R = A^{(0)}A^{(1)}\cdots A^{(N-1)}$ and $M = A^{(N)} A^{(N+1)} \cdots A^{(N+p-1)}$.
	As a shorthand, we shall use the following notation $\vec{v} = R\overline{M}$.

	Below, when we mention a MCF algorithm, we mean a MCF algorithm for some given $(\mathcal{I},\A)$ and $n$.

	\begin{example}
Let $\mathcal{I} = \{\R^3_+\}$, $\A = \{A\}$ with $A = \begin{pmatrix}
1&1&0\\
0&1&0\\
0&0&1
\end{pmatrix}$. In this $(\mathcal{I},\mathcal{A})$ $2$-dimensional continued fraction algorithm, every vector has an expansion equal to $(A,A,\dots,)$. 
However, no expansion is periodic since we never have $\vec{v}^{(N)} = \lambda \vec{v}^{(N+p)}$  for some $N \geq 0, p >0$ and $\lambda \in \R$.
	\end{example}

As we have seen from the example above, a careful choice of the sets $\mathcal{I}$ and $\A$ is crucial for the properties of the algorithm. Therefore, we work mainly with algorithms in which every expansion corresponds to a uniquely given vector (up to a multiple). Examples of these algorithms follow later in this section.

	\subsection{Transvections} \label{sec:transvections}

	Let $\SL{}(n,\Z)$ be the special linear group of matrices over $\Z$ of dimension $n \times n$ with determinant~$1$.
	Let $\SL{}(n,\N)$ be the subset of $\SL{}(n,\Z)$ containing all the matrices with nonnegative elements.
In what follows, we focus mainly on $\A \subseteq \SL{}(n,\N)$.
	However, the monoid $\SL{}(n,\N)$ is not finitely generated for $n\geq 3$ (for a proof, see Chapter 12.5 of \cite{fogg2002substitutions}).
	On the other hand, the group $\SL{}(n,\Z)$ is finitely generated by \emph{transvections} which are matrices $T_{ij}$ that have $1$'s on the diagonal and on the $i,j$-th position and $0$'s elsewhere. 
	It follows that for $n = 3$ we have
	\[
	T_{12} = \begin{pmatrix}
		1 & 1 & 0 \\
		0 & 1 & 0 \\
		0 & 0 & 1
	\end{pmatrix},\quad
	T_{13} = \begin{pmatrix}
		1 & 0 & 1 \\
		0 & 1 & 0 \\
		0 & 0 & 1
	\end{pmatrix},\quad
	T_{21} = \begin{pmatrix}
		1 & 0 & 0 \\
		1 & 1 & 0 \\
		0 & 0 & 1
	\end{pmatrix},
	\]
	\[
	T_{23} = \begin{pmatrix}
		1 & 0 & 0 \\
		0 & 1 & 1 \\
		0 & 0 & 1
	\end{pmatrix},\quad
	T_{31} = \begin{pmatrix}
		1 & 0 & 0 \\
		0 & 1 & 0 \\
		1 & 0 & 1
	\end{pmatrix},\quad
	T_{32} = \begin{pmatrix}
		1 & 0 & 0 \\
		0 & 1 & 0 \\
		0 & 1 & 1
	\end{pmatrix}.
	\]
	
	The following result due to Conder, Robertson and Williams (\cite{conder1992presentations}) gives us even the presentation of $\SL{}(n,\Z)$.
	Let $[A,B]$ be the \emph{commutator of $A$ and $B$}, i. e., $[A,B] = ABA^{-1}B^{-1}$.

	\begin{proposition}[\cite{conder1992presentations}] \label{prop:gen_SL}
		The group $\SL{}(n,\Z)$, where $n\geq 3$, has a presentation with the $n(n-1)$ generators $T_{ij}$ subject only to the Steinberg relations
		\begin{gather}
		\begin{aligned}
		\left[ T_{ij},T_{jk} \right] &= T_{ik} && \text{for } i \neq k,\\
		\left[ T_{ij},T_{k\ell} \right] &= 1          && \text{for } i \neq \ell, j \neq k,\\
		&&& \text{where } i,j,k \in \{1,\dots,n\}\\
		\text{and to the relation } (T_{12}T_{21}^{-1}T_{12})^4 &= 1.
		\end{aligned}
		\end{gather}
	\end{proposition}

	\subsection{Regular continued fractions}
	
	We present several well-known Jacobi--Perron type multidimensional continued fraction algorithms and some of their properties.
	The lowest dimension that we can study is $n-1 = 1$. 
	In this dimension, we usually consider only one algorithm and that is the regular continued fraction algorithm.

	We set $\mathcal{I}_{RCF} = \{I_{RCF,1},I_{RCF,2}\}$, $\A_{RCF} = \{C_1,C_2\}$
	with $I_{RCF,1} = \left\{ \vek{v_1}{v_0} \colon v_1 \geq v_0 > 0 \right\}$, $I_{RCF,2} = \left\{ \vek{v_1}{v_0} \colon 0< v_1 < v_0 \right\}$ and
	\[
	C_1 = \mat{1}{1}{0}{1} \quad \text{ and } \quad  C_2 = \mat{1}{0}{1}{1}.
	\]
	
	\begin{definition} \label{def:regular_CF}
		\emph{Regular continued fractions} are the $(\mathcal{I}_{RCF},\A_{RCF})$ one-dimensional continued fractions.
	\end{definition}
	
	Now, we describe the concrete form of \Cref{alg:MCF_alg} for the sets $\mathcal{I}_{RCF}$ and $\A_{RCF}$ in more detail. In this case, the $i$-th step of the algorithm for $\vec{v} = \begin{pmatrix}
		v_1\\
		v_0
	\end{pmatrix} \in \R_{+}^{2}$ is as follows:

	\begin{enumerate}
		\item if $ \vec{v}^{(i)} \in I_{RCF,1}$, then we set:
		\[
		\vec{v}^{(i+1)} \coloneqq C_1^{-1} \vec{v}^{(i)} = \vek{v_1^{(i)} - v_0^{(i)}}{v_0^{(i)}},
		\]
		\[
		C^{(i)} \coloneqq C_1;
		\]
		\item if $ \vec{v}^{(i)} \in I_{RCF,2}$, then we set:
		\[
		\vec{v}^{(i+1)} \coloneqq C_2^{-1}  \vec{v^{(i)}} = \vek{v_1^{(i)}}{v_0^{(i)} - v_1^{(i)}},
		\]
		\[
		C^{(i)} \coloneqq C_2.
		\]
	\end{enumerate}

	This is the so-called additive form of regular continued fraction algorithm.
    This algorithm is very often presented in its multiplicative form:
	
	\begin{definition} \label{def:regular_CF_multiplicative}
		\emph{Multiplicative regular continued fractions} are the $(\mathcal{I}_{RCFM},\A_{RCFM})$ one-dimensional continued fractions with
		$\mathcal{I}_{RCFM} = \{I_{RCFM,1,k},I_{RCFM,2,k}\colon k \in \Z_+\}$, where for all $k \in \Z_+$ we have 
		$I_{RCFM,1,k} = \left\{ \vek{v_1}{v_0}  \in \R_{+}^{2} \colon \left\lfloor \frac{v_1}{v_0}\right\rfloor = k \right\}$ and
		$I_{RCFM,2,k} = \left\{ \vek{v_1}{v_0}  \in \R_{+}^{2} \colon \left\lfloor \frac{v_0}{v_1}\right\rfloor = k \right\}$, and
		$\A_{RCFM} = \{C_{RCFM,1,k},C_{RCFM,2,k}\colon k \in \Z_+\}$ with $C_{RCFM,1,k} = \mat{1}{k}{0}{1} = C_1^k$ and $C_{RCFM,´2,k} = \mat{1}{0}{k}{1} = C_2^k$.
		
	\end{definition}

	\begin{remark}
	For a given $\vek{v_1}{v_0} \in \R_{+}^{2}$, using the additive form of regular continued fraction algorithm, we obtain the expansion
	\[
	\vek{v_1}{v_0} = \left( C^{(0)}, C^{(1)}, \dots  \right) = \left( \underbrace{C_1,\dots,C_1}_{a_0 \text{ times}}, \underbrace{C_2,\dots,C_2}_{a_1 \text{ times}}, \dots \right)
	\]
	with $a_0 \in \N$ and  $a_\ell \in \Z_+$ for $\ell > 0$.
	This correspond to the multiplicative regular continued fraction expansion as follows:
	$\vek{v_1}{v_0} = \left( C_{1}^{a_0}, C_{2}^{a_1}, \dots \right)$ if $a_0 > 0$ and $\vek{v_1}{v_0} = \left(  C_{2}^{a_1},  C_{2}^{a_2} \dots \right)$ if $a_0 = 0$.
	The reason to fix the matrix $C_1$ as the first is to obtain the correspondence with the following form of regular continued fraction representation:
	\[ \frac{v_1}{v_0} = a_0+ \cfrac{1}{a_1+\cfrac{1}{a_2+\cfrac{1}{\ddots}}}, \]
which is also usually written as $[a_0,a_1,a_2,\dots]$.
	\end{remark}
	
	In what follows, we sum up the properties of regular continued fraction algorithm that we will investigate also in the multidimensional case.	
	As already mentioned, regular continued fraction expansions can be used to determine some algebraic properties of the components of the represented vector $\vek{v_1}{v_0}$.
	We can determine whether the number $\frac{v_1}{v_0}$ has algebraic degree equal to $1$, $2$ or greater than $2$.
	For $\frac{v_1}{v_0} \in \Q$, the regular continued fraction expansion of $\vek{v_1}{v_0}$ is finite.
	For $\frac{v_1}{v_0} \in \R_+ \setminus \Q$, the algorithm does not terminate, and the regular continued fraction expansion of the vector $\vek{v_1}{v_0}$ is unique.
	Moreover, the regular continued fraction expansion of a vector $\vek{v_1}{v_0}$ is eventually periodic if and only if $\frac{v_1}{v_0}$ is a quadratic irrational number.
	It is purely periodic if and only if $\frac{v_1}{v_0}$ is a quadratic irrational, $\frac{v_1}{v_0} >1$ and for the algebraic conjugate $\left(\frac{v_1}{v_0}\right)'$ of $\frac{v_1}{v_0}$ we have 
	$\left(\frac{v_1}{v_0}\right)' \in (-1,0)$.

	\begin{example}
	We have
		\[
		\vek{1}{2} = \left( C_2,C_1\right)
		\]
		in the multiplicative regular continued fraction algorithm.

		Since
		\[
		\vek{\sqrt{2}}{1} = C_1 \vek{\sqrt{2}-1}{1} = C_1(C_2^2 C_1^2)^{m}  \vek{\sqrt{2}-1}{1} \text{ for all }  m \in \N,
		\]
		we have
		\[
		\vek{\sqrt{2}}{1} = (C_1,\overline{C_2^2, C_1^2}) \quad \text{ and } \quad \vek{\sqrt{2} - 1}{1} = (\overline{C_2^2, C_1^2})
		\]
		in the multiplicative regular continued fraction algorithm.
\end{example}
	
	Now, we show how the situation for $n-1\geq 2$ differs from $n-1 = 1$ and describe the commonly used Jacobi--Perron type MCF algorithms and their basic properties.

	\subsection{The Jacobi--Perron algorithm}
	\label{sec:JP_alg}
	The earliest and probably most extensively studied MCF algorithm is the Jacobi--Perron algorithm (the JPA). It was introduced by Jacobi in \cite{jacobi} (1868\footnote{Published by E. Heine from the legacy of G. G. J. Jacobi. Jacobi studied the question of MCF at least since 1839 (and died in 1851).}) for dimension $n-1 = 2$. And Perron (1907 \cite{Perron}) generalised it to an arbitrary dimension.

	The algorithm works only with vectors $\vec{v} = \begin{pmatrix}
		v_1\\
		\vdots\\
		v_{n}\\
	\end{pmatrix}\in \R_+^{n}$, where $v_n > v_\ell$ for all $\ell \in \{1,\dots,n-1\}$. We put
	\[
	\mathcal{I}_{JP} = \{I_{JP,j_2,\dots,j_{n}} \subseteq \R_+^{n}\colon  j_\ell\in \N \text{ for } \ell \in \{2,\dots,n-1\}\text{ and } j_n \in \Z_{+}\},
	\]
	where
	\[
	\vec{v} \in I_{JP,j_2,\dots,j_{n}} 
	\iff \left \lfloor \frac{v_\ell}{v_1} \right \rfloor = j_\ell \text{ for all } \ell \in \{2,\dots, n\}
	\]
	and
	\[
	\A_{JP} = \{A_{JP,j_2,\dots,j_{n}} \subseteq \Z^{n,n}\colon  j_\ell\in \N \text{ for } \ell \in \{2,\dots,n-1\}\text{ and } j_n \in \Z_{+}\},
	\]
	where
	\[\quad A_{JP,j_2,\dots,j_{n}} = \prod_{\ell = 2}^{n}
	T_{1\ell}^{j_{\ell}}P,
	\]
	the matrices $T_{1\ell}$ are transvections (defined in \Cref{sec:transvections}) and $P$ is a permutation matrix such that
	\[
	P_{\ell m} = \begin{cases}
		1 & \text{ iff } \left(\ell \in \{2,\dots,n\} \wedge m = \ell -1 \right ) \vee \left(\ell =1 \wedge m = n \right) \\
		0 & \text{ otherwise }
	\end{cases}
	\]
	(i. e. it shifts the first column to the end).
	
	Specially for $n = 3$, we obtain
	\[
	A_{JP,j_2,j_3} = \begin{pmatrix}
		0 & 0 & 1\\
		1 & 0 & j_2\\
		0 & 1 & j_3
	\end{pmatrix}.
	\]
	This means that in the $i$-th step of the algorithm the vector $\vec{v}^{(i)}$ is transformed as follows:
	\	\[
	(v_1^{(i)},v_2^{(i)},v_3^{(i)}) \mapsto \left ( v_2^{(i)} - \left \lfloor \frac{v_2^{(i)}}{v_1^{(i)}} \right \rfloor v_1^{(i)}, v_3^{(i)}- \left \lfloor \frac{v_3^{(i)}}{v_1^{(i)}} \right \rfloor v_1^{(i)},v_1^{(i)} \right) = (v_1^{(i+1)}, v_2^{(i+1)}, v_3^{(i+1)}).
	\]
	\\
	
	It is easy to verify that $\mathcal{I}_{JP}$ is a set of disjoint subsets of $\R_+^n´$ and that for every $\alpha >0$ and every $I_i \in \mathcal{I}_{JP}$ we have $\alpha I_i \subseteq I_i$. 
	
	\begin{definition}[Jacobi (1868 \cite{jacobi}) and Perron (1907 \cite{Perron})] \label{def:JPA}
		
		The \emph{Jacobi--Perron algorithm (or simply the JPA)} is the $(\mathcal{I}_{JP},\A_{JP})$ MCF algorithm.
		
		The \emph{Jacobi--Perron expansion (or simply the JP expansion)} of a vector is the $(\mathcal{I}_{JP},\A_{JP})$ MCF expansion.
	\end{definition}

\begin{example}
	In the Jacobi--Perron algorithm, we have:
	\[
\begin{pmatrix}
	1\\
	\sqrt[3]{2}\\
	\sqrt[3]{4}\\
\end{pmatrix} = A_{JP,1,1}A_{JP,2,3}\overline{A_{JP,3,3}}.
\]
\end{example}

	In \cite{Perron} Perron proved that every sequence of matrices from $\A_{JP}$ represents a unique vector (up to scalar multiplication).
	Moreover, David \cite{david1957contribution} showed that for $n-1 = 2$ the JPA detects rational dependence. And Brentjes \cite{brentjes} showed an example for $n-1 \geq 3$ in which the JPA fails to detect rational dependence.

	\subsection{The Brun algorithm}
	Brun (1920 \cite{brun1920}) proposed another Jacobi--Perron type algorithm.
	
	The algorithm works with vectors $\vec{v} = \begin{pmatrix}
		v_1\\
		\vdots\\
		v_n
	\end{pmatrix} \in \R_+^{n}$. We have
	\[
	\mathcal{I}_{B} =\{I_{\sigma}\colon \sigma \in S_n\} 
	\quad \text{with}\quad
	I_{\sigma} = \left \{\begin{pmatrix}
		v_1\\
		\vdots\\
		v_n
	\end{pmatrix} \in \R_+^{n} \colon 0<v_{\sigma(1)}<v_{\sigma(2)}<\dots<v_{\sigma(n)} \right\}
	\]

	and
	\[
	\mathcal{A}_{B} = \{A_{B\sigma}\subseteq \Z^{n,n}\colon  \sigma \in S_n\} \quad \text{ with }\quad
	A_{B\sigma} = T_{\sigma(n)\sigma(n-1)},
	\]
	where the matrices $T_{jk}$ ($j,k \in \{1,\dots,n\}$) are transvections (defined in \Cref{sec:transvections}).

	Especially for $n-1 = 2$ we obtain that $
	\A_{B} = \{A_{B(123)},A_{B(132)},\dots,A_{B(321)}\}$, where $
	A_{B(321)}=
	T_{12},
	A_{B(231)}=T_{13},
	A_{B(312)} = 
	T_{21},
	A_{B(132)}= T_{23},
	A_{B(213)} = 
	T_{31}\text{ and } A_{B(123)} = 
	T_{32}.$
	Therefore, the $i$-th step of the algorithm works as follows. If $\vec{v}^{(i)} \in I_{(123)}$, then
	\[
	(v_1^{(i)},v_2^{(i)},v_3^{(i)})^T \mapsto \left (v_1^{(i)}, v_2^{(i)}, v_3^{(i)} - v_2^{(i)} \right)^T = (v_1^{(i+1)}, v_2^{(i+1)}, v_3^{(i+1)}),
	\]
	and analogously for $\vec{v}^{(i)}$ in other sets from $\mathcal{I}_{B}$.
	
	We can easily see that $\mathcal{I}_B$ is a set of disjoint subsets of $R_+^n$ and that for every $\alpha >0$ and every $I_i \in \mathcal{I}_B$ we have $\alpha I_i \subseteq I_i$. 
	
	\begin{definition}[Brun (1920 \cite{brun1920})] \label{def:Brun}
		The \emph{Brun algorithm} is the $(\mathcal{I}_{B},\A_{B})$ MCF algorithm.
		
		The \emph{Brun expansion} of a vector is the $(\mathcal{I}_{B},\A_{B})$ MCF expansion.
	\end{definition}
	
	\begin{example}
		In the Brun algorithm, we have:
		\[
		\begin{pmatrix}
			1\\
			\sqrt[3]{2}\\
			\sqrt[3]{4}\\
		\end{pmatrix} = T_{32}\overline{ T_{21}T_{13}^{3}T_{32}T_{23}^{3}T_{32}T_{21}^{3}T_{13}T_{31}T_{12}T_{23}T_{31}T_{12}}.
		\]
	\end{example}

	Brun in \cite{brun1920} proved that for the dimension $n-1=2$ every expansion in the Brun algorithm corresponds to a uniquely determined vector (up to scalar multiplication). This result was later generalised by Greiter \cite{greiter1977mehrdimensionale} for arbitrary dimension.
	
	For dimension $n-1 = 2$, Brun provided another important result and that is that the Brun algorithm detects rational dependence for $n-1 = 2$. (Analogue of this theorem probably does not hold for any other $n -1 \geq 3$.)

	\subsection{The Selmer algorithm}
	In response to the Brun algorithm, Selmer in 1961 \cite{selmer} proposed the following Jacobi--Perron type algorithm.
	
	The algorithm works with vectors $\vec{v} = \begin{pmatrix}
		v_1\\
		\vdots\\
		v_n
	\end{pmatrix} \in \R_+^{n}$. We have
	\[
	\mathcal{I}_{S}= \{I_{\sigma}\colon \sigma \in S_n\} 
	\quad \text{with}\quad
	I_{\sigma} = \left \{\begin{pmatrix}
		v_1\\
		\vdots\\
		v_n
	\end{pmatrix} \in \R_+^{n} \colon 0<v_{\sigma(1)}<v_{\sigma(2)}<\dots<v_{\sigma(n)} \right\}
	\]

	and
	\[
	\mathcal{A}_{S} = \{A_{S\sigma}\subseteq \Z^{n,n}\colon \sigma \in S_n\} \quad \text{ with }\quad
	A_{S\sigma} = T_{\sigma(n)\sigma(1)},
	\]
	where the matrices $T_{jk}$ ($j,k \in \{1,\dots,n\}$) are transvections (defined in \Cref{sec:transvections}).

Especially for $n-1 = 2$, we obtain that $
\A_{S} = \{A_{S(123)},A_{S(132)},\dots,A_{S(321)}\}$, where $
A_{S(123)}=
T_{31},
A_{S(132)}=T_{21},
A_{S(213)} = 
T_{32},
A_{S(231)}= T_{12},
A_{S(312)} = 
T_{23}\text{ and } A_{S(321)} = 
T_{13}.$

Therefore, the $i$-th step of the algorithm works as follows. If $\vec{v}^{(i)} \in I_{(123)}$, then
\[
(v_1^{(i)},v_2^{(i)},v_3^{(i)})^T \mapsto \left (v_1^{(i)}, v_2^{(i)}, v_3^{(i)} - v_1^{(i)},  \right)^T = (v_1^{(i+1)}, v_2^{(i+1)}, v_3^{(i+1)}),
\]
and analogously for $\vec{v}^{(i)}$ in other sets from $\mathcal{I}_{S}$.

\begin{definition}[Selmer (1961 \cite{selmer})]
	The \emph{Selmer algorithm} is the $(I_{S},\A_{S})$ MCF algorithm.
	
	The \emph{Selmer expansion} of a vector is the $(I_{S},\A_{S})$ MCF expansion.
\end{definition}

\begin{example}
	In the Selmer algorithm, we have:
	\[
	\begin{pmatrix}
		1\\
		\sqrt[3]{2}\\
		\sqrt[3]{4}\\
	\end{pmatrix} = T_{31}\overline{T_{23}T_{13}T_{21}T_{32}T_{12}T_{31}T_{21}T_{32}T_{13}T_{23}T_{12}T_{32}T_{13}T_{21}T_{31}}.
	\]
\end{example}

\section{Matrices of multiplication in a number field} \label{sec:wekla_conv}
We focus on matrices of multiplication $T_\lambda$ in a number field $\Q(\alpha)$ (as a vector space over $\Q$) of degree $n$.
In this section, we show that the transpose of such a matrix is fully determined by any of its single columns.
Moreover, we show that the mappings which determine the matrix from a single selected column are linear.

We start with a lemma on eigenvalues of the transposed matrix $T_{\lambda}^{\vec{v}}$.

\begin{lemma} \label{le:eigenvalue_mult_transpose}
		Let $\vec{v} = \begin{pmatrix}
			v_1\\
			\vdots\\
			v_n
		\end{pmatrix}$ be a basis of some number field, $\lambda \in \Q(v_1,\dots,v_n)$ and $T_\lambda^{\vec{v}}$ be the  matrix of the linear transformation $t_{\lambda}$ in the basis $\vec{v}$. We have
	
	\begin{equation*} \label{eq:transposition}
	M \vec{v} = \lambda \vec{v} \iff M = (T_\lambda^{\vec{v}})^T.
	\end{equation*}
\end{lemma}
\begin{proof}
Let $\vec{e_i}$ be the $i$-th vector of the standard basis, i.e., $(e_i)_j = \begin{cases} 1 & \text {if } i =j; \\ 0 & \text{otherwise}.\end{cases}$. It follows from the definition of the matrix $T_\lambda^{\vec{v}}$ that
\[
\vec{v}^T T_\lambda^{\vec{v}} \vec{e_i} = \lambda v_i
\]
and therefore
\[
 \vec{e_i}^T (T_\lambda^{\vec{v}})^T \vec{v} = \lambda v_i.
\]
This holds for every $i \in \{1,\dots,n\}$ and therefore
\[
(T_\lambda^{\vec{v}})^T \vec{v} = \lambda \vec{v}.\qedhere
\]
\end{proof}

Since we are later (in \Cref{sec:periodic_exp}) interested in this exact situation of $\lambda$ being an eigenvalue of a matrix $M$, we state the next theorem with the transposition of $T_{\lambda}^{\vec{v}}$.

	\begin{theorem} \label{th:weak_exists_Q}
	Let  $\vec{v} = \begin{pmatrix}
		v_1 \\
		\vdots \\
		v_{n}
	\end{pmatrix}$ be a basis (of a finite field extension of degree $n$ as a vector space over $\Q$), $\ell \in \left\{ 1,\dots,n \right\}$ and $\lambda \in \Q(v_1,\dots,v_{n})$.
	There exists a mapping $\mathcal{Q}_{\ell, \vec{v}} : \R^n \mapsto \R^{n,n}$ such that
	for every $\lambda \in \Q(v_1,\dots,v_n)$
	we have
	\[
	M = \mathcal{Q}_{\ell,\vec{v}
	} \left( M_{\bullet,\ell}  \right)
	\]
	where ${M} = {T_{\lambda}^{\vec{v}}}^T$.

	Moreover, there exists an $n$-tuple $Q_{\ell,\vec{v}}$ of matrices from $\Q^{n,n}$ such that
	their $i$-th component satisfies
	\[
	\left( Q_{\ell,\vec{v}} \right)_i  M_{\bullet,\ell} =
	\left(\mathcal{Q}_{\ell,\vec{v}} \left( M_{\bullet,\ell}  \right)\right)_{\bullet,i}.
	\]
\end{theorem}

\begin{proof} 	
	By~\Cref{le:eigenvalue_mult_transpose}, we have that $\lambda$ is an eigenvalue of $M$ corresponding to the eigenvector $\vec{v}$.
	For all $i \in \{1,\dots,n\}$, we have 
	\begin{equation} \label{eq:Mv}
	M_{i,\ell}  = ((\lambda v_i)_{\vec{v}})_\ell
	\end{equation}
	where $(\lambda v_i)_{\vec{v}}$ denotes the vector of coordinates of $\lambda v_i$ in the basis $\vec{v}$.
	We show by contradiction that the eigenvalue $\lambda$ is uniquely determined by these equations.

	Let $\lambda_1$ and $\lambda_2$ be distinct numbers for which \eqref{eq:Mv} holds.
	Set $\lambda_3 = \lambda_1 - \lambda_2 \neq 0$.
	We have
	\begin{equation} \label{eq:v_ell_zero}
	0 = ((\lambda_3 v_i)_{\vec{v}})_\ell
	\end{equation}
	for all $i \in \{1,\dots,n\}$.
	
	Since $\lambda_3 \neq 0$, the linear transformation $t_{\lambda_3}$ is an automorphism of $\Q(v_1,\dots,v_n)$, hence its matrix is regular.
	On the other hand, equality \eqref{eq:v_ell_zero} implies that its matrix in the basis $\vec{v}$ has zeros on the $\ell$-th row, which is a contradiction.

	Hence $\lambda$ is uniquely determined by \eqref{eq:Mv}, i.e., it can be determined from $\vec{v}$ and $M_{\bullet,\ell}$.
	Therefore, we can also find the whole matrix $M = {T_{\lambda}^{\vec{v}}}^T$.	
	
	The moreover part follows from the fact that the elements of $M$ are linear combinations of the coordinates of $\lambda$ in the basis $\vec{v}$.
	\end{proof}

	\begin{remark}
		Let $\alpha \in \R \setminus \left\{ 0 \right\} $. If $M$ is a matrix of a linear transformation $t_{\lambda}$ in the basis $\vec{v}$ (as a vector space over $\Q$), then it is also a matrix of the same linear transformation in the basis $\alpha \vec{v}$. 
		
		Therefore, we have
		\[
		\mathcal{Q}_{\ell,\vec{v}} = \mathcal{Q}_{\ell,\alpha \vec{v}}
		\]
		for all $\alpha \in \R \setminus \left\{ 0 \right\} $.
	\end{remark}

	In what follows, we keep the same notation as in \Cref{th:weak_exists_Q}, i.e., we assign to the mapping $\mathcal{Q}_{\ell,\vec{v}
		}$ the $n$-tuple of matrices $Q_{\ell,\vec{v}}$.
		We demonstrate this and the claim of \Cref{th:weak_exists_Q} in the following example.

	\begin{example}
		Let $\vec{v} = \tvek{\sqrt[3]{4}}{\sqrt[3]{2}}{1}$.
		We have 
		\[
		Q_{1,\vec{v}} =\left( \begin{pmatrix}
		1&0&0\\
		0&1&0\\
		0&0&1\\
		\end{pmatrix},
		\begin{pmatrix}
		0&0&2\\
		1&0&0\\
		0&1&0\\
		\end{pmatrix},
		\begin{pmatrix}
		0&2&0\\
		0&0&2\\
		1&0&0
		\end{pmatrix} \right).
	\]
		
		It follows that
		\[
		\mathcal{Q}_{1,\vec{v}} \left ( \tvek{x}{y}{z} \right ) = \begin{pmatrix}
		x&2z&2y \\
		y&x&2z\\
		z&y&x
		\end{pmatrix} 
		\]

		We take three matrices: $M= \begin{pmatrix}
		1&2&2\\
		1&1&2\\
		1&1&1
		\end{pmatrix}$, $M^2$ and $M^3$. The matrix $M$ is a transposition of a matrix of linear transformation $t_\varepsilon$ in the basis $\vec{v}$, where $\varepsilon = (\sqrt[3]{4}+ \sqrt[3]{2}+1)$ is a unit in $\Ocal_{\Q(\sqrt[3]{2})}$. It follows that $M^2$ and $M^3$ are transpositions of the matrix of transformation $t_{\varepsilon^2}$ and $t_{\varepsilon^3}$ respectively.
		
		We have
		\begin{align*}
		M &=  \begin{pmatrix}
		1&2&2\\
		1&1&2\\
		1&1&1
		\end{pmatrix}
		= \mathcal{Q}_{1,\vec{v}} \left ( M_{\bullet,1} \right ) =\mathcal{Q}_{1,\vec{v}} \left ( \tvek{1}{1}{1} \right ) \\
		&= \left( (Q_{1,\vec{v}})_1 M_{\bullet,1}\quad (Q_{1,\vec{v}})_2 M_{\bullet,1}\quad (Q_{1,\vec{v}})_3 M_{\bullet,1} \right)
		\\ &=
		 \begin{pmatrix}
		1+ 0 + 0&0+0+2&0+2+0\\
		0+1+ 0&1+0+0&0+0+2\\
		0+0+1&0+1+0&1+0+0
		\end{pmatrix},
		\\
		M^2 &=  \begin{pmatrix}
		5&6&8\\
		4&5&6\\
		3&4&5
		\end{pmatrix}
		=\mathcal{Q}_{1,\vec{v}} \left ( \tvek{5}{4}{3} \right ) \\ 
		&=  \left( (Q_{1,\vec{v}})_1 M^2_{\bullet,1}\quad (Q_{1,\vec{v}})_2 M^2_{\bullet,1}\quad (Q_{1,\vec{v}})_3 M^2_{\bullet,1} \right) \\
		&=  \begin{pmatrix}
		5+ 0 + 0&0+0+2 \cdot 3&0+2 \cdot 4+0\\
		0+4+ 0&1 \cdot 5+0+0&0+0+2 \cdot 3\\
		0+0+3&0+1 \cdot 4+0&1 \cdot 5+0+0
		\end{pmatrix},
		\\
		M^3 &=  \begin{pmatrix}
		19&24&30\\
		15&19&24\\
		12&15&19
		\end{pmatrix}
		=\mathcal{Q}_{1,\vec{v}} \left ( \tvek{19}{15}{12} \right ) \\
		&= \left( (Q_{1,\vec{v}})_1 M^3_{\bullet,1} \quad (Q_{1,\vec{v}})_2 M^3_{\bullet,1}\quad (Q_{1,\vec{v}})_3 M^3_{\bullet,1} \right) \\
		&=  \begin{pmatrix}
		19+ 0 + 0&0+0+2 \cdot 12&0+2 \cdot 15+0\\
		0+15+ 0&1 \cdot 19+0+0&0+0+2 \cdot 12\\
		0+0+12&0+1 \cdot 15+0&1 \cdot 19+0+0
		\end{pmatrix}.
		\end{align*}
	\end{example}

The next lemma shows that the mapping $\mathcal{Q}_{\ell, \vec{v}}$ of \Cref{th:weak_exists_Q} can be determined from the first $n$ powers of a matrix of multiplication.

\begin{lemma} \label{lem:Qi_from_M_repetend}
	Let $\vec{v} = \begin{pmatrix}
		v_1 \\
		\vdots\\
		v_{n}\\
	\end{pmatrix}$ be a basis of $\Q(v_1)$ as a vector space over $\Q$ and $\lambda \in \Q(v_1)$ be an algebraic number of degree $n$.
	
	We can determine the elements of the matrices of the $n$-tuples $Q_{\ell,\vec{v}} $ (for $\ell \in \{1,\dots,n\}$) as linear combinations of elements of the first $n$ powers of the matrix $M = {T_\lambda^{\vec{v}}}^T$.
\end{lemma}
\begin{proof}
	Firstly, we realise that $\lambda$ is an eigenvalue of $M$ and $\vec{v}$ is the corresponding eigenvector. The degree of $\lambda$ is $n$ and therefore $M$ has $n$ distinct eigenvalues.
	
	The elements of the matrices of $n$-tuples $Q_{\ell,\vec{v}} $ are linear combinations of the minimal polynomial of $v_1$, the coordinates of $v_2,\dots, v_n$ in the basis $1,v_1,\dots,v_1^{n-1}$ and the coordinates of $\lambda$ in the basis $\vec{v}$. This means that the $n$-tuples $Q_{\ell,\vec{v}} $ are uniquely determined by the matrix $M$.
	
	Now we show that we can obtain the elements of the matrices  of these $n$-tuples as linear combinations of the elements in the $\ell$-th columns of the first $n$ powers of $M$.
	
	We need to find all the $n^3$ elements of the matrices of the $n$-tuple $Q_{\ell,\vec{v}}$. These elements are given by the $n^3$ equations that we obtain by expressing the elements of the matrices $M, M^2,\dots, M^n$ as linear combinations of the elements in their $\ell$-th columns. This follows directly from the definition of the $n$-tuple $Q_{\ell,\vec{v}}$.
	
	We show by contradiction that this system of linear equations is nonsingular. Suppose otherwise. The existence of the $n$-tuple $Q_{\ell,\vec{v}}$ implies that there exist at least two solutions of this system of linear equations. Using the definition of the $n$-tuples $Q_{\ell,\vec{v}}$, we obtain that there is a vector $\vec{x} \in \Q^{n}$, $\vec{x} \neq 0$, such that $\vec{x}^T (M^m)_{\bullet,\ell}  = 0$ for all $m \in \{1,\dots,n\}$. This means that we have an equation with $n$ variables and $n$ solutions. It implies that the  solutions $ (M)_{\bullet,\ell},\dots, (M^n)_{\bullet,\ell}$ are linearly dependent.
	Now, let $\vec{w} = (w_1,\dots,w_n)$ be a left eigenvector of $M$ corresponding to an eigenvalue $\beta$. 
	It follows that also $\beta$ has degree $n$ and that $\vec{w}$ is a basis of $\Q(v_1)$. The elements of $(M^m)_{\bullet,\ell}$ are in fact coordinates of $\beta^m w_\ell$ in the basis $\vec{w}$. This means that $\beta w_\ell,\dots,\beta^n w_\ell$ are linearly dependent. At the same time, $\beta$ has degree $n$ and therefore $w_\ell = 0$. This is a contradiction.
\end{proof}

	The next theorem refines \Cref{th:weak_exists_Q} to the case when the components of $\vec{v} = \begin{pmatrix}
		y^{n-1}\\
		\vdots\\
		y\\
		1
		\end{pmatrix}$ form a polynomial basis.
		In this case, we can explicitly determine the matrices of $\mathcal{Q}_{\ell,\vec{v}}$ by the coefficients of the monic minimal polynomial of $y$.
		For simplicity, we state this claim for a specific value of $\ell$, namely for $\ell = 1$.
		For other values of $\ell$, analogous formulas can be obtained.

	\begin{theorem} \label{prop:T_epsilon}
		Let $y$ be an algebraic number of degree $n$ such that
		$
		\sum_{r=0}^{n-1} \alpha_r y^r + y^{n} = 0,
		$
		where $\alpha_r \in \Q$, and	
		$\vec{v} = \begin{pmatrix}
		y^{n-1}\\
		\vdots\\
		y\\
		1
		\end{pmatrix}$.

		Let $i,j,k \in \{1,\dots n\}$. We have
		\[
		\left(\left( Q_{1,\vec{v}} \right)_i \right)_{j,k} = 
		\begin{cases} 1 & \text{ for } i \leq j , k = j-i+1 \\
		\alpha_{n-i+1+j-k} & \text{ for } 2 \leq i \leq j, k \in \{j-i+2,\dots, j\} \\
		- \alpha_{n-i+1+j-k} & \text{ for } j<i , j+1\leq k \leq n+j-i+1 \\
		0 & \text{ otherwise }
		\end{cases}.
		\]

	\end{theorem}
	\begin{proof}
		Suppose that we have a matrix ${T_\lambda^{\vec{v}}}^T$.
			Because $\vec{v}$ is a basis of a finite field extension and $\lambda \in \Q(y)$, we can find numbers $\beta_1,\dots,\beta_n \in \Z$ such that
		 $\lambda = \sum_{j = 0}^{n-1}\beta_jy^j$. We put $\beta_i = 0$ for all $i <0$.
		
		From the definition of the matrix $T_{\lambda}^{\vec{v}}$ we obtain
		\[
		\left( T_\lambda^{\vec{v}} \right)_{i,j} =
		\displaystyle \beta_{-i+j} + \sum_{k = j}^{n-1}\beta_k  \sum_{r =0}^{\min \{ k-j, n-i\}}
		\alpha_{n-j-r} \sum_{\mathclap{\substack{p_1\dots p_m  \\ m \leq k-j - r\\ p_s \geq 1\\ \sum_{s=1}^m p_s = k-j - r}}} (-1)^{m+1} \alpha_{n-p_1}\dots \alpha_{n-p_m}.
		\]
		
		Moreover, if we put $x_j = \left( T_\lambda^{\vec{v}} \right)_{1,j}$ and $\alpha_j = 0$ for all $j <0$ we obtain
		
		\begin{equation} \label{eq:T_exponential_12}
		\left( T_\lambda^{\vec{v}} \right)_{1,j} = 
		x_j, \quad
		\left( T_\lambda^{\vec{v}} \right)_{2,j} = x_{j-1} + x_j \alpha_{n-1}
		\end{equation}
		
		for all $j \in \{1,\dots, n\}$,
		\begin{equation} \label{eq:T_exponential_i_leq_j}
		\left( T_\lambda^{\vec{v}} \right)_{i,j} = x_{j-i + 1} + x_{j-i+2} \alpha_{n-1} + \dots + x_j \alpha_{n-i+1}
		\end{equation}
		for all $i \in \{3,\dots,n\} ,j \in \{1,\dots,n\}, i \leq j$ and
		\begin{equation}\label{eq:T_exponential_i_g_j}
		\left( T_\lambda^{\vec{v}} \right)_{i,j} = \sum_{m = 1}^{n-j} - x_{j+m} \alpha_{n-i-m+1}
		\end{equation}
		for all $i,j \in \{1,\dots,n\}, i > j$. 
		
		Now it remains to realize that $\left(\left( Q_{1,\vec{v}} \right)_i \right)_{j,k} $ is equal to the coefficient of $x_k$ in~\eqref{eq:T_exponential_i_g_j}.
		\end{proof}

\subsection{Refinements on $Q_{\ell,\vec{v}}$ for $n=3$}
\label{sec:extracting_information}

The explicit formulas for the elements of $\left( Q_{\ell,\vec{v}} \right)_i$ can be derived without relying on the specific form of $\vec{v}$ in the previous theorem. 

However, as the approach of the previous proof would yield more intricate expressions, we focus in this subsection on the case $n=3$ for which we state more general claims.

\begin{lemma}\label{prop:T_epsilon_cubic_case}
	Let $y$ be a cubic number such that $\alpha_0 + \alpha_1 y + \alpha_2 y^2 + y^3 = 0$, where $\alpha_0, \alpha_1,\alpha_2 \in \Q$, $x = \gamma_0 + \gamma_1 y + \gamma_2 y^2$, where $\gamma_0, \gamma_1, \gamma_2 \in \Q$ and $\vec{v} = \tvek{x}{y}{1} $.
	
	We have $Q_{1,\vec{v}} = \left(\begin{pmatrix}
		1&0&0\\
		0&1&0\\
		0&0&1\\
	\end{pmatrix},
	\begin{pmatrix}
		0 & b_1 & c_1 \\
		1 & b_2 & c_2 \\
		0 & b_3 & c_3 \\
	\end{pmatrix},
	\begin{pmatrix}
		0 & c_1 & c_4 \\
		0 & c_2 & c_5 \\
		1 & c_3 & c_6
	\end{pmatrix} \right) $, where

	\begin{equation} \label{eq:b1_c_1_minpoly_new}
	\begin{aligned}
		\gamma_0 &= - c_2, & \gamma_1 &= c_3-b_2, & \gamma_2 &= b_3,\\
		\\
		\alpha_0 &= \frac{-c_1-c_3c_2}{b_3^2}, & \alpha_1 &= \frac{c_3^2-b_1-b_3c_2-b_2c_3}{b_3^2}, &  \alpha_2 &= \frac{2c_3 - b_2}{b_3}.
	\end{aligned}	
	\end{equation}
	and
	\begin{equation} \label{eq:c_4_new}
		c_4  = \frac{c_1c_3-c_1b_2+c_2b_1}{b_3}, \
		c_5 = \frac{c_3c_2+c_1}{b_3} \ \text{ and } \ c_6 = \frac{c_3^2+b_3c_2-b_1-b_2c_3}{b_3}.
	\end{equation}
	
	Or equivalently
	\begin{equation} \label{eq:b_1_vyjadrene_minpoly_new}
		\begin{aligned}
			b_1 &= \gamma_2\gamma_0+\gamma_1\alpha_2\gamma_2-\gamma_1^2-\alpha_1\gamma_2^2
			\\
			b_2 &= \alpha_2\gamma_2-2\gamma_1
			\\
			b_3 &= \gamma_2
			\\
			c_1 &= \gamma_0\alpha_2\gamma_2-\gamma_0\gamma_1-\alpha_0\gamma_2^2
			\\
			c_2 &=-\gamma_0
			\\
			c_3 &= \alpha_2\gamma_2-\gamma_1
			\\
			c_4 &=\gamma_0\alpha_1\gamma_2-\gamma_0^2-\alpha_0\gamma_2\gamma_1
			\\
			c_5 &=-\alpha_0\gamma_2
			\\
			c_6 &=\alpha_1\gamma_2-2\gamma_0.
		\end{aligned}
	\end{equation}		
\end{lemma}

If we choose to normalize some other component of $\vec{v}$, the mapping $\mathcal{Q}$ is transformed by a suitable permutation, which yields the following claim:

\begin{proposition} \label{prop:permutations}
	Let $y$ be a cubic number for which we have $\alpha_0 + \alpha_1 y + \alpha_2 y^2 + y^3 = 0$, where $\alpha_0, \alpha
	_1,\alpha_2 \in \Q$, and $x = \gamma_0 + \gamma_1 y + \gamma_2 y^2$, where $\gamma_0, \gamma_1, \gamma_2 \in \Q$, $\vec{v_2} = \tvek{1}{x}{y}$, resp. $\vec{v_3} = \tvek{y}{1}{x}$.

	We obtain that 
	\[ Q_{2,\vec{v_2}} = \left( \begin{pmatrix}
		c_6 & 1 & c_3  \\
		c_4& 0 & c_1  \\
		c_5& 0 & c_2 \\
	\end{pmatrix},
	\begin{pmatrix}
		1&0&0\\
		0&1&0\\
		0&0&1\\
	\end{pmatrix},
	\begin{pmatrix}
		c_3& 0 & b_3 \\
		c_1& 0 & b_1 \\
		c_2& 1 & b_2\textsl{} \\
	\end{pmatrix} \right),
	\]
	\[ \text{resp.} \quad
	Q_{3,\vec{v_3}} = \left(\begin{pmatrix}
		b_2 & c_2 & 1 \\
		b_3 & c_3 & 0\\
		b_1 & c_1 & 0\\
	\end{pmatrix},
	\begin{pmatrix}
		c_2 & c_5 & 0 \\
		c_3 & c_6 & 1\\
		c_1 & c_4 & 0\\
	\end{pmatrix},
	\begin{pmatrix}
		1&0&0\\
		0&1&0\\
		0&0&1\\
	\end{pmatrix} \right),\]

	where \eqref{eq:b1_c_1_minpoly_new}, \eqref{eq:c_4_new} and \eqref{eq:b_1_vyjadrene_minpoly_new} hold.
\end{proposition}

\begin{proof}
Let $P = \begin{pmatrix}
0 & 0 & 1 \\ 1 & 0 & 0 \\ 0 & 1 & 0
\end{pmatrix}$ be the permutation matrix determined by $\vec{v_2} = P \vec{v}$ where $\vec{v} = \tvek{x}{y}{1}$ as in~\Cref{prop:T_epsilon_cubic_case}.
Let $\pi$ be the permutation given by the permutation matrix $P$; we have $\pi(1) = 2$.
It follows from the definition of $\mathcal{Q}_{\ell,\vec{v}}$ and \Cref{prop:T_epsilon_cubic_case} that
\begin{equation} \label{eq:per_Q}
Q_{\pi(1),\vec{v_2}} = \left(  P \left(  Q_{1,\vec{v}} \right)_1 P^{-1}, P \left(  Q_{1,\vec{v}} \right)_2 P^{-1}, P \left(  Q_{1,\vec{v}} \right)_3 P^{-1}     \right)  P^T. 
\end{equation}
Applying $P$ yields the first part of the desired result.

The second part is obtained analogously: we have $\vec{v_3} = \widetilde{P} \vec{v}$ for $\widetilde{P} = \begin{pmatrix}
0 & 1 & 0 \\ 0 & 0 & 1 \\ 1 & 0 & 0
\end{pmatrix}$.
\end{proof}

Now we show that we can obtain some non-trivial information about the vector $\vec{v}$ directly from the $3$-tuples $Q_{\ell,\vec{v}}$.

\begin{proposition} \label{prop:comparison_of_elements}
	Let $\vec{v} = \tvek{x}{y}{1} \in \R_+^3$ be a basis of some complex cubic number field (as a vector space over $\Q$).
	Moreover, let
	\[Q_{1,\vec{v}} = \left(\begin{pmatrix}
		1&0&0\\
		0&1&0\\
		0&0&1\\
	\end{pmatrix},
	\begin{pmatrix}
		0 & b_1 & c_1 \\
		1 & b_2 & c_2 \\
		0 & b_3 & c_3 \\
	\end{pmatrix},
	\begin{pmatrix}
		0 & c_1 & c_4 \\
		0 & c_2 & c_5 \\
		1 & c_3 & c_6
	\end{pmatrix} \right),\]
	\[ Q_{2,\vec{v}} = \left( \begin{pmatrix}
		\widetilde{c_6} & 1 & \widetilde{c_3} \\
		\widetilde{c_4} & 0 & \widetilde{c_1} \\
		\widetilde{c_5} & 0 & \widetilde{c_2} \\
	\end{pmatrix},
	\begin{pmatrix}
		1&0&0\\
		0&1&0\\
		0&0&1\\
	\end{pmatrix},
	\begin{pmatrix}
		\widetilde{c_3} & 0 & \widetilde{b_3}\\
		\widetilde{c_1} & 0 & \widetilde{b_1} \\
		\widetilde{c_2} & 1 & \widetilde{b_2} \\
	\end{pmatrix} \; \right)
	\] and  
	\[ \; Q_{3,\vec{v}} = \left( \begin{pmatrix}
		\widehat{b_2} & \widehat{c_2} & 1 \\
		\widehat{b_3} & \widehat{c_3} & 0 \\
		\widehat{b_1} & \widehat{c_1} & 0 \\
	\end{pmatrix},
	\begin{pmatrix}
		\widehat{c_2} & \widehat{c_5} & 0 \\
		\widehat{c_3} & \widehat{c_6} & 1 \\
		\widehat{c_1} & \widehat{c_4} & 0 \\
	\end{pmatrix},
	\begin{pmatrix}
		1&0&0\\
		0&1&0\\
		0&0&1\\
	\end{pmatrix} \right).
	\]
	
	We have
	\[
	b_3(b_3+2c_3-b_2-c_5+c_6-2c_2)> 0 \iff y<1,
	\]
	\[
	\widetilde{b_3}(\widetilde{b_3}+2\widetilde{c_3}-\widetilde{b_2}-\widetilde{c_5}+\widetilde{c_6}-2\widetilde{c_2})> 0 \iff x>1
	\]
	and
	\[
	\widehat{b_3}(\widehat{b_3}+2\widehat{c_3}-\widehat{b_2}-\widehat{c_5}+\widehat{c_6}-2\widehat{c_2})> 0 \iff x<y.
	\]
\end{proposition}

\begin{proof}
Let $f$ be the monic minimal polynomial of $y$, i.e.,
	\[
	f(y) = y^3 + \alpha_2y^2+\alpha_1y+\alpha_0 = 0
	\]
	for $\alpha_2,\alpha_1,\alpha_0 \in \Q$.
As $y$ is the only real root of the polynomial $f$, the function $f$ is strictly increasing.
Hence,
	\[
	y < 1 \iff 0 = f(y) < f(1) = 1 + \alpha_2 + \alpha_1 + \alpha_0.
	\]
	Using \eqref{eq:b_1_vyjadrene_minpoly_new}, this is equivalent to
	\[
	0 < 1 + \frac{2c_3 - b_2}{b_3} + \frac{c_6 - 2 c_2}{b_3} - \frac{c_5}{b_3},
	\]
	which proves the first part of the statement.
	Proofs of the other two parts of the statement are analogous; they rely on \Cref{prop:permutations} and  \Cref{rem:renormalisation}.
\end{proof}

The last claim requires the cubic number field to be complex.
We indicate here that this is indeed necessary.

	Let $\vec{v} = \tvek{x}{y}{1} \in \R_+^3$ be a basis of some cubic number field (as a vector space over $\Q$).
	We set $\alpha_2,\alpha_1,\alpha_0,\gamma_2,\gamma_1,\gamma_0 \in \Q$, $\gamma_2 \neq 0$ to be such that
	\begin{equation} \notag
		y^3+\alpha_2y^2+\alpha_1y +\alpha_0 = 0
	\end{equation}
	and
	\begin{equation} \notag
		x = \gamma_2y^2 +\gamma_1y+\gamma_0.
	\end{equation}
	
	It can happen that there is more than one positive real vector $\vec{v}$ for which these two equalities hold. 
	In such a case, the mapping $\mathcal{Q}_{1,\vec{v}}$ is not an injection.

	If $y_1,y_2$ are real roots of $y^3 + \alpha_2y^2+\alpha_1y+\alpha_0$ for which there exist $x_1,x_2 \in \R^{+}$ such that $ x_1 = \gamma_2y_1^2 + \gamma_1 y_1 + \gamma_0, x_2 = \gamma_2y_2^2 + \gamma_1y_2 + \gamma_0$
	and $0<y_1<1, 1<y_2$, then we cannot decide whether $y>1$ only from the knowledge of the triplet $Q_{1,\vec{v}}$ (without the knowledge of $\vec{v}$).
	
	The situation for other inequalities is analogous.

\begin{proposition} \label{prop:norms_of_elements}
	Let $\vec{v} = \tvek{x}{y}{1} \in \R^3_+$ be a basis of some cubic number field (as a vector space over $\Q$).
	Moreover, let
	
	\[Q_{1,\vec{v}} = \left(\begin{pmatrix}
		1&0&0\\
		0&1&0\\
		0&0&1\\
	\end{pmatrix},
	\begin{pmatrix}
		0 & b_1 & c_1 \\
		1 & b_2 & c_2 \\
		0 & b_3 & c_3 \\
	\end{pmatrix},
	\begin{pmatrix}
		0 & c_1 & c_4 \\
		0 & c_2 & c_5 \\
		1 & c_3 & c_6
	\end{pmatrix} \right),\]
	\[ Q_{2,\vec{v}} = \left( \begin{pmatrix}
		\widetilde{c_6} & 1 & \widetilde{c_3} \\
		\widetilde{c_4} & 0 & \widetilde{c_1} \\
		\widetilde{c_5} & 0 & \widetilde{c_2} \\
	\end{pmatrix},
	\begin{pmatrix}
		1&0&0\\
		0&1&0\\
		0&0&1\\
	\end{pmatrix},
	\begin{pmatrix}
		\widetilde{c_3} & 0 & \widetilde{b_3}\\
		\widetilde{c_1} & 0 & \widetilde{b_1} \\
		\widetilde{c_2} & 1 & \widetilde{b_2} \\
	\end{pmatrix} \; \right)
	\] and  
	\[ \; Q_{3,\vec{v}} = \left( \begin{pmatrix}
		\widehat{b_2} & \widehat{c_2} & 1 \\
		\widehat{b_3} & \widehat{c_3} & 0 \\
		\widehat{b_1} & \widehat{c_1} & 0 \\
	\end{pmatrix},
	\begin{pmatrix}
		\widehat{c_2} & \widehat{c_5} & 0 \\
		\widehat{c_3} & \widehat{c_6} & 1 \\
		\widehat{c_1} & \widehat{c_4} & 0 \\
	\end{pmatrix},
	\begin{pmatrix}
		1&0&0\\
		0&1&0\\
		0&0&1\\
	\end{pmatrix} \right).
	\]
	
	We have
	\begin{equation}
		\left | \frac{N(y)}{N(1)} \right | = \left | \frac{c_5}{b_3}\right| , \quad
		\left | \frac{N(1)}{N(x)} \right | = \left | \frac{\widetilde{c_5}}{\widetilde{b_3}}\right|, \quad \text{ and } 
		\left | \frac{N(x)}{N(y)} \right | = \left | \frac{\widehat{c_5}}{\widehat{b_3}}\right|. 
	\end{equation}
	
\end{proposition}

\begin{proof}
	The numbers $x,y,1$ constitute a basis of a cubic number field, and therefore, $y$ is a cubic number. We put
	$
	y^3 + \alpha_2y^2+\alpha_1y+\alpha_0 = 0.
	$
	We know that \eqref{eq:b_1_vyjadrene_minpoly_new} holds. Therefore,
	\[
	\left | \frac{N(y)}{N(1)} \right | = |\alpha_0| =  \left | \frac{c_5}{b_3}\right|.
	\]
	
	The rest of the proof is analogous using \Cref{prop:permutations}, \Cref{rem:renormalisation}, the triplet $Q_{2,\vec{v}}$ (resp. $Q_{3,\vec{v}}$) and the minimal polynomial of $\widetilde{y} = \frac{1}{x} $ (resp. $\widehat{y} =\frac{x}{y} $).
\end{proof}

\section{Periodic MCF expansions} \label{sec:periodic_exp}

In this section, we consider eventually periodic expansions of $\vec{v}$, i.e., we assume $\vec{v} = R\overline{N}$ (where both the preperiodic part and the repetend are represented as a product).
We shall not distinguish between purely and eventually periodic expansions by considering the matrix $RNR^{-1}$, called the \emph{matrix of repetend}, and the equality $R\overline{N} = \overline{RNR^{-1}}$.
If $\vec{v} = \overline{M}$, i.e.,  the matrix $M$ is a matrix of a repetend of an expansion of $\vec{v}$, then $\vec{v}$ is an eigenvector of $M$.
We use the following reformulation of this fact from \cite{brentjes}:

\begin{theorem}[\cite{brentjes}, Theorem 3.1.]
	\label{thm:periodicity_of_MCF}
	Let $\vec{v} =\begin{pmatrix}
		v_1\\
		\vdots\\
		v_{n}
	\end{pmatrix} \in \R_+^{n}$, $\vec{v} = \overline{M}$ in a given unimodular MCF algorithm. We have
	\[
	M \vec{v} = \lambda  \vec{v},
	\]
	where $\lambda \in \R$. Moreover,
	\begin{itemize}
		\item $\lambda$ is an algebraic unit of degree at most $n$;
		\item If the degree of $\lambda$ equals $n$, then the numbers $\frac{v_1}{v_n},\dots,\frac{v_{n-1}}{v_n},\frac{v_n}{v_n}$ constitute a basis (as a vector space over $\Q$) of the number field $\Q(\lambda)$. 
	\end{itemize}
\end{theorem}

We cannot omit the condition on the degree of $\lambda$ since $\mathrm{deg}( \lambda) \leq n-1$ would allow $\frac{v_j}{v_n} \not \in \Q(\lambda)$. For an example of such a vector and algorithm, see Remark (1) in \cite{brentjes}.

The following theorem states that the matrix of repetend always equals to a matrix of multiplication by some unit in basis $\vec{v}$.

	\begin{theorem} \label{thm:candidates_onMP_algebraic_new}
		Let $\vec{v} = \begin{pmatrix}
		y_{1}\\
		\vdots \\
		y_{n}\\
		\end{pmatrix}$ be a basis of $\Q(y_1)$ (as a vector space over $\Q$), where $\vec{v}$ has an eventually periodic expansion in a unimodular Jacobi--Perron type MCF algorithm.
		Moreover, let $M$ be a matrix of repetend of this MCF expansion of $\vec{v}$. We have
		\[
		M = {T^{\vec{v}}_\varepsilon}^T,
		\]
		where $\varepsilon \in U(\mathcal O_{\Q(y_1)})$ and $T^{\vec{v}}_\varepsilon $ is a matrix of linear transformation $t_{\varepsilon}$ (defined by \eqref{eq:t_alpha}) in the basis $\vec{v}$.
	\end{theorem}
	\begin{proof}
		It follows from \Cref{thm:periodicity_of_MCF} that
		\begin{equation} \label{eq:Mv_epv}
		M \vec{v} = \varepsilon \vec{v},
		\end{equation}
		where $\varepsilon$ is an algebraic unit. 
		Moreover, the matrix $M$ is an integer matrix, and therefore, for all $i \in \{1,\dots,n\}$ we have $(M\vec{v})_i \in \Q(y_1)$, hence $\varepsilon \in \Q(y_1)$.
		Equality $M =	{T^{\vec{v}}_\varepsilon}^T$ follows from~\Cref{le:eigenvalue_mult_transpose}.
	\end{proof}

	\begin{remark} \label{rem:relations_in_T_epsilon}
		Let $\vec{v} = \begin{pmatrix}
		v_1\\
		\vdots\\
		v_{n}\\
		\end{pmatrix}$ be a basis of a number field of degree $n$ (as a vector space over $\Q$), $\lambda, \widehat{\lambda} \in \Q(v_1,\dots,v_n) $ and $m \in \Z$. We have
		\begin{equation}
		\label{eq:T_epsilon_commute}
		T_\lambda^{\vec{v}} T_{\widehat{\lambda}}^{\vec{v}} = T_{\widehat{\lambda}}^{\vec{v}} T_\lambda^{\vec{v}} \; \text{ and } \; T_{\lambda^m}^{\vec{v}} = \left(T_\lambda^{\vec{v}}\right)^m,
		\end{equation}
		which implies that $\{ {T_\lambda^{\vec{v}}}^T| \lambda \in \Q(v_1,\dots,v_n), \lambda \neq 0\}$ is an Abelian group.
		\end{remark}

\subsection{Purely periodic MCF expansions}

In this subsection, we focus on purely periodic Jacobi--Perron type MCF expansions and state two necessary conditions for a MCF expansion to be purely periodic.
First of all, we introduce the weak convergence of the Jacobi--Perron type MCF algorithms.

\begin{definition}
	Let $\left(  M^{(s)} \right)_{s=0}^{+\infty}$
	be a sequence of matrices from $\R^{n,n}$.
	Moreover, let $j \in \{1,\dots,n\}$.
	We say it \emph{weakly converges to $\vec{v} \in \R^{n}$ with respect to the $j$-th column} if the following two conditions are fulfilled:
	\begin{enumerate}
		\item there exists $\widetilde{P}$ such that $M^{(P)}$ is positive for all $P > \widetilde{P}$;
		
		\item the sequence
		\[
		\left( \frac{M^{(s)}_{i,j}}{M^{(s)}_{k,j}} \right)_{s=P}^{+\infty} 
		\]
		converges to $\frac{\vec{v}_i}{\vec{v}_k}$ for all $i \in \{1,\dots,n\}$ and some $k \in \{1,\dots,n\}$.
		
	\end{enumerate}
\end{definition}

\begin{remark}
	Since all elements of all matrices $M^{(s)}$ for $s \geq \widetilde{P}$ are positive, we can choose the integer $k$ arbitrarily.
\end{remark}

The $(\mathcal{I},\A)$ $(n-1)$-dimensional MCF algorithm is \emph{weakly convergent} if for every vector $\vec{v} \in \R_+^{n}$ 
whose expansion is $ \left(  A^{(0)}, A^{(1)}, \dots \right)$ with $M^{(s)} = A^{(0)}A^{(1)} \cdots A^{(s)}$ we have that the sequence $M^{(s)}$ weakly converges to $\vec{v}$ with respect to the $j$-th column for every $j$.

\begin{remark}
Note that the definition of weak convergence varies in the literature.
This is mostly due to variances in the definitions of MCF algorithms themselves.
We base our definition of weak convergence on the definition in the book \cite{brentjes} of Brentjes although the definition present there is based on the geometric definition of MCF algorithm.
The reader may also refer to Schweiger~\cite{schweiger2000} on the matter of various concepts of convergence.	
\end{remark}

Examples of Jacobi--Perron type algorithms that are weakly convergent in arbitrary dimension are the Jacobi--Perron algorithm (proved by Perron in \cite{Perron}) and the Brun algorithm (proved by Greiter in \cite{greiter1977mehrdimensionale}).

Let us now recall that a matrix $M$ is \emph{primitive} if there exists a positive integer $k$ such that every element of $M^k$ is positive.
In this article, we use the Perron--Frobenius theorem, which gives us key information about the eigenvectors of a primitive matrix. 
We state the theorem in a form suitable for the rest of the article:

\begin{theorem}[Perron--Frobenius theorem] \label{thm:PF}
	Let $M$ be a primitive matrix.
	\begin{enumerate}
		\item The matrix $M$ has a positive real eigenvalue $\lambda_{\max}$ such that every other eigenvalue $\lambda$ satisfies
		\[
		|\lambda| < \lambda_{\max}.
		\]
		\item The eigenvalue
		$\lambda_{\max}$ has algebraic and geometric multiplicity equal to one and has an eigenvector $\vec{v}$ such that every component of $\vec{v}$ is positive.
		\item Any eigenvector with nonnegative components is a multiple of $\vec{v}$.
	\end{enumerate}
\end{theorem}

In the following, we will discuss vectors composed of conjugates of components of $\vec{v}$ and expansions of such vectors.
More precisely, consider a vector $\vec{v}^T=({v_1},\dots, {v_n})$ and let $\alpha$ be a primitive element for the number field generated by $\vec{v}$, i.e., $\Q({v_1},\dots, {v_n})=\Q(\alpha)$. Fixing a conjugate $\widehat{\alpha}$ of $\alpha$ determines an embedding $\Q(\alpha)\hookrightarrow \C$ given by $\beta\mapsto\widehat\beta$ where $\widehat \beta = g(\widehat \alpha)$ with $\beta = g(\alpha)$ for some polynomial $g$ (over $\Q$).
We set $\widehat{\vec{v}}^T = (\widehat{v_1},\dots, \widehat{v_n})$.
It follows that
\[
{T_{\varepsilon}^{\vec{v}}} = {T_{\widehat{\varepsilon}}^{\widehat{\vec{v}}}}.
\]

\begin{theorem} \label{prop:purly_per}
	Let $\vec{v} \in \R^n_+$ be a basis of a number field of degree $n$ (as a vector space over $\Q$) and $\widehat{\vec{v}} \neq \vec{v}$ be a conjugate vector of $\vec{v}$.
	Moreover, suppose that $\vec{v}$ has a \emph{purely} periodic expansion in some unimodular weakly-convergent $(n-1)$-dimensional continued fraction algorithm.
	We have
	\[
	\widehat{\vec{v}} \not \in \R^n_+.
	\]
\end{theorem}
\begin{proof}

	Let $\vec{v} = \overline{M}$.
	Moreover, as the algorithm is weakly convergent, the matrix $M$ is primitive.
	
	By~\Cref{thm:candidates_onMP_algebraic_new} there exists a unit $\varepsilon \in U(\mathcal O_{\Q(\vec{v})}) $ such that $M = {T_{\varepsilon}^{\vec{v}}}^T = {T_{\widehat{\varepsilon}}^{\widehat{\vec{v}}}}^T$.
	It follows that both $\vec{v}$ and $\widehat{\vec{v}}$ are distinct eigenvectors of $M$.
	Moreover $\widehat{\vec{v}}$ is not a multiple of $\vec{v}$ since that would imply $\varepsilon = \widehat \varepsilon$, which is impossible as $\widehat{\vec{v}} \neq \vec{v}$.
	As $M$ is primitive and $\vec{v} \in \R^n_+$, by the Perron--Frobenius~\Cref{thm:PF} any eigenvector with nonnegative components is a multiple of $\vec{v}$.
	Therefore, $\widehat{\vec{v}} \not \in \R^n_+$.

\end{proof}

	\begin{proposition} \label{cor:not_purely_per}
		Suppose that $\vec{v} =  \begin{pmatrix}
			y^{n-1}\\
			\vdots\\
			y\\
			1
		\end{pmatrix}$, where $y$ is an algebraic number of degree $n$.
		If $\vec{v}$ has a \emph{purely} periodic expansion in some unimodular $(\mathcal{I},\mathcal{A})$  $(n-1)$-dimensional continued fraction algorithm for which $\mathcal{A}\subset \SL{}(n,\N)$, then the norm $N(y)$ has sign $(-1)^{n-1}$.

	\end{proposition}
	\begin{proof}
Let $\Id$ be the identity matrix of dimension $n$ and $\alpha_0,\dots,\alpha_{n-1} \in \Q$ be such that 
\[
\sum_{j=0}^{n-1} \alpha_j y^j + y^{n} = 0.
\]

We suppose for contradiction that the sign of $N(y)$ is $(-1)^{n}$. This is equivalent to saying that $\alpha_0=(-1)^nN(y) >0$.

We show that $\{ T_\varepsilon^{\vec{v}}|\varepsilon \in U(\Ocal_{\Q(y)})\} \cap \mathrm{SL}(n,\N) = \left\{ \Id \right\}$ (and therefore also $\{ {T_\varepsilon^{\vec{v}}}^T|\varepsilon \in U(\Ocal_{\Q(y)})\} \cap \mathrm{SL}(n,\N) = \left\{ \Id \right\}$).
A purely periodic expansion in a MCF algorithm has a matrix of repetend that is equal to a product of matrices from the set $\mathcal{A}$ and by \Cref{thm:candidates_onMP_algebraic_new}, the matrix of repetend is equal to the matrix $T_{{\varepsilon}}^{\vec{v}}$ for some ${\varepsilon} \in U(\Ocal_{\Q(y)})$.	
	Since $\mathcal{A} \subset \SL{}(n,\N)$, we have $T_{\varepsilon}^{\vec{v}} \in \mathrm{SL}(n,\N)$.
By \Cref{eq:T_exponential_i_g_j} we obtain
		\[
		(T_\varepsilon^{\vec{v}})_{n,j-1} = - (T_\varepsilon^{\vec{v}})_{1,j} \alpha_0,
		\]
		for all $j \in \{2,\dots, n\}$.
Since $\alpha_0 > 0$ and $T_{\varepsilon}^{\vec{v}} \in \mathrm{SL}(n,\N)$, it follows that $(T_\varepsilon^{\vec{v}})_{1,j} = 0$ for all $j \in \{2,\dots,n\}$.
Using \Cref{eq:T_exponential_12} and \Cref{eq:T_exponential_i_leq_j}, we obtain that $T_{\varepsilon}^{\vec{v}} = \Id$.
Since $\mathcal{A} \subset \SL{}(n,\N)$ and no product of matrices from $\A$ is equal to $\Id$ (by \Cref{def:MCF}), no matrix of repetend is equal to $\Id$, thus we have a contradiction.
	\end{proof}

	\section{Candidates on the matrix of repetend}
	\label{sec:candidates}

	In this section,	
	we first show how to generate all the matrices $T_\varepsilon^{\vec{v}}$ for every $\vec{v} = \begin{pmatrix}
v_1\\
\vdots\\
v_{n-1}\\
1	
\end{pmatrix}$ and every $\varepsilon \in U(\mathcal O_{\Q(v_{n-1})})$.

The procedure we describe relies on~\Cref{lem:Qi_from_M_repetend} and requires knowledge of the minimal polynomial of $v_{n-1}$, the coordinates of the components of $\vec{v}$ in the polynomial basis $(1,v_{n-1},\dots,v_{n-1}^{n-1})$ (of $\Q(v_{n-1})$), and the fundamental units of $\mathcal O_{\Q(v_{n-1})}$.
\Cref{thm:candidates_onMP_algebraic_new} shows that every matrix of repetend of $\vec{v}$ in some unimodular Jacobi--Perron type MCF algorithm can be expressed as ${T_\varepsilon^{\vec{v}}}^T$ for some $\varepsilon \in U(\mathcal O_{\Q(v_{n-1})})$.
Therefore, we refer to the matrices ${T_\varepsilon^{\vec{v}}}^T$ as the \emph{candidates on the matrix of repetend}.
Again, we limit our exposition to the case of $n=3$ for simplicity but note that the procedure generalizes to larger values of $n$.

\subsection{Finding candidates on the matrix of repetend}
\label{sec:finding_candidates}

	Let $\widehat{y}$ be a cubic number for which $\alpha_0 + \alpha_1 \widehat{y} + \alpha_2 \widehat{y}^2 + \widehat{y}^3 = 0$, where $\alpha_0, \alpha_1,\alpha_2 \in \Q$, $\widehat{x} = \gamma_0 + \gamma_1 \widehat{y} + \gamma_2 \widehat{y}^2$, where $\gamma_0, \gamma_1, \gamma_2 \in \Q$ and $\vec{v} = \tvek{\widehat{x}}{\widehat{y}}{1}$ be a basis of some cubic number field (as a vector space over $\Q$). 
	We continue by the description of a procedure finding all the candidates on the matrix of repetend of the MCF expansion of the vector $\vec{v}$.

	Firstly, we have to realise that the number $ \widehat{y}$ is a cubic number, and therefore, by Dirichlet's \Cref{thm:dirichlet}, there are either one or two fundamental units in $\mathcal O_{\Q( \widehat{y})}$.

	Let 
\begin{equation} \label{eq:varepsilon_in_basis}
	\varepsilon_1 = \beta_1+\beta_2\widehat{y} + \beta_3\widehat{x}, \text{ resp. } \varepsilon_1 = \beta_1+\beta_2\widehat{y} + \beta_3\widehat{x}, \varepsilon_2 = \widehat{\beta}_1+\widehat{\beta}_2\widehat{y} + \widehat{\beta}_3\widehat{x},
\end{equation} be the fundamental unit, resp. units, of $\mathcal O_{\Q( \widehat{y})}$.

It follows from \Cref{thm:candidates_onMP_algebraic_new,thm:roots_of_unity} and \Cref{eq:unit,eq:T_epsilon_commute} that every 
	candidate $M$ on the matrix of repetend of the MCF expansion of $\tvek{\widehat{x}}{\widehat{y}}{1}$ can be written as
	\begin{equation} \label{eq:Mc}
	M = \pm \left ({T_{\varepsilon_1}^{\vec{v}}}^T\right)^{m_1} 
	\end{equation}
	for $m_1 \in \Z$, respectively
	\begin{equation} \label{eq:Mr}
	M = \pm \left ( \left ({T_{\varepsilon_1}^{\vec{v}}}^T \right)^{m_1}\left({T_{\varepsilon_2}^{\vec{v}}}^T\right)^{m_2} \right )
	\end{equation}
	for $m_1,m_2 \in \Z$.

	We can easily verify by direct computation that $\left({T_{\varepsilon_1}^{\vec{v}}}^T\right)_{\bullet,1} =
	\tvek{x_1}{y_1}{z_1}$, where 
	\begin{align} \label{eq:eps1_to_first_column}
		x_1 &= \beta_1+ \beta_2\left(\frac{\gamma_1}{\gamma_2} -\alpha_2\right) + \beta_3\left(\frac{\gamma_1^2}{\gamma_2}-2\gamma_1\alpha_2 + \alpha_2^2\gamma_2+2\gamma_0-\alpha_1\gamma_2\right), \nonumber \\ y_1 &= \frac{\beta_2}{\gamma_2} + \beta_3\left(\frac{\gamma_1}{\gamma_2}-\alpha_2\right), \\
		z_1 &= \beta_3 \nonumber
	\end{align}
	and, applicable in the case of two fundamental units, $\left({T_{\varepsilon_2}^{\vec{v}}}^T\right)_{\bullet,1} = \tvek{x_2}{y_2}{z_2}$ where 
	\begin{align} \label{eq:eps2_to_first_column}
		x_2 &= \widehat{\beta}_1+ \widehat{\beta}_2\left(\frac{\gamma_1}{\gamma_2} -\alpha_2\right) + \widehat{\beta}_3\left(\frac{\gamma_1^2}{\gamma_2}-2\gamma_1\alpha_2 + \alpha_2^2\gamma_2+2\gamma_0-\alpha_1\gamma_2\right), \nonumber \\
		y_2 &= \frac{\widehat{\beta}_2}{\gamma_2} + \widehat{\beta}_3\left(\frac{\gamma_1}{\gamma_2}-\alpha_2\right),\\
		z_2 &= \widehat{\beta}_3. \nonumber
	\end{align}

	Now, we can use \Cref{th:weak_exists_Q} and \Cref{prop:T_epsilon_cubic_case} to compute the matrices ${T_{\varepsilon_1}^{\vec{v}}}^T$, resp. ${T_{\varepsilon_1}^{\vec{v}}}^T$ and ${T_{\varepsilon_2}^{\vec{v}}}^T$.
	We use the notation from \Cref{prop:T_epsilon_cubic_case}. We obtain that 
	\[
	{T_{\varepsilon_1}^{\vec{v}}}^T = \left( \left(Q_{1,\vec{v}} \right)_1 \tvek{x_1}{y_1}{z_1} \quad \left(Q_{1,\vec{v}} \right)_2 \tvek{x_1}{y_1}{z_1} \quad \left(Q_{1,\vec{v}} \right)_3 \tvek{x_1}{y_1}{z_1} \right)
	\]
	and similarly for the matrix ${T_{\varepsilon_2}^{\vec{v}}}^T$.
	
	For simplicity, we do the explicit calculation only for the case $\widehat{x} = \widehat{y}^2$. In this case, we obtain a simpler form, and that is $x_1 = \beta_1-\beta_3\alpha_1-\beta_2\alpha_2+\beta_3\alpha_2^2, y_1 = \beta_2-\beta_3\alpha_2, z_1 = \beta_3$ and eventually $x_2 = \widehat{\beta}_1-\widehat{\beta}_3\alpha_1-\widehat{\beta}_2\alpha_2+\widehat{\beta}_3\alpha_2^2, y_2 = \widehat{\beta}_2-\widehat{\beta}_3\alpha_2, z_2 = \widehat{\beta}_3$.
	
	For $i \in \{1,2\}$, we obtain that 
	\begin{equation} \label{eq:prvni_sloupec}
	{T_{\varepsilon_i}^{\vec{v}}}^T = \begin{pmatrix}
	x_i & -\alpha_1y_i-\alpha_0z_i &-\alpha_0y_i\\
	y_i & x_i+\alpha_2y_i&-\alpha_0z_i\\
	z_i & y_i+\alpha_2z_i&x_i+\alpha_2y_i+\alpha_1z_i
	\end{pmatrix}.
	\end{equation}
	Note that the last equality holds also if $\varepsilon_i$ is not a unit.

In the most common case, when $\mathcal{A} \subseteq \SL{}(n,\N)$, some of the candidate matrices $M$ (given by \eqref{eq:Mc}, resp. \eqref{eq:Mr}) can be excluded.
First, the determinant of the matrix of repetend has to be $1$, hence we can exclude $M$ if it has determinant equal to $-1$.
Second, the matrix of repetend has integer entries.
Thus, we can exclude $M$ if it has non-integer entries.
If $\widehat{y}$ is an algebraic integer, then $M$ has always integer entries.

It remains to comment on the fact that the knowledge of fundamental units is required to find all the candidates on the matrix of repetend.

The procedure of obtaining all fundamental units of a real quadratic number field is known (for example, see \cite[Section 11]{alaca2004introductory}).
If $\alpha = \sqrt[3]{d}$ for some $d \in \N$, $d \neq e^3$ where $e \in \N$, the problem of finding fundamental units in $\Ocal_{\Q(\alpha)}$ is closely connected with the cubic analogue of the Pell's equation, and therefore, we can use the process described in \cite{barbeau2006pell}.
For some other algebraic number fields, there are algorithms for computing a set of fundamental units. Most of these algorithms are based on the geometric interpretation of MCFs. One of these algorithms is the Voronoi's algorithm for computing a set of fundamental units of a cubic number field (1896, \cite{voronoi1896generalization} and later restated in a different form in \cite{delone1964theory}). In 1985, Buchmann (\cite{buchmann1985generalization} and \cite{buchmann1985generalizationII}) generalized Voronoi's algorithm to an arbitrary number field with the group of units of rank $1$ and $2$.
For an example of sets of fundamental units in some cubic number fields see \cite[Section 13.6]{alaca2004introductory} or \cite{cusick1987table}.

We illustrate the described procedure of finding all candidates on a matrix of repetend on the following example:

\begin{example} \label{ex:matrix_of_repetend}
	Let $\widehat{y}$ be the only positive root of the polynomial $\widehat{y}^3 + \widehat{y}^2-2 \widehat{y} - 1$. We investigate the MCF expansion of the vector $\vec{v} = \tvek{\widehat{y}^2}{\widehat{y}}{1}$. Since the number $\widehat{y}$ has three real conjugates, by~\Cref{thm:dirichlet} there are two fundamental units in $U(\mathcal O_{\Q(\widehat{y}, \widehat{y}^2)})$. In \cite{alaca2004introductory} we can find that the two fundamental units are $\varepsilon_1 = -1 + \widehat{y} + \widehat{y}^2$ and $\varepsilon_2 = 2-\widehat{y}^2$.
	
	Using the computation in \eqref{eq:Mr}, we obtain that every candidate $M$ on the matrix of repetend is defined by
	\[
	M = \pm (M_1^{m_1} M_2^{m_2})
	\]
	where $m_1,m_2 \in \Z$,
	\[
	M_1 = \begin{pmatrix}
	1&1&0\\
	0&1&1\\
	1&1&-1
	\end{pmatrix} \quad \text{and} \quad
	M_2 = \begin{pmatrix}
	-1&1&1\\
	1&0&-1\\
	-1&0&2
	\end{pmatrix}.
	\]

We compare it with the expansion of $\vec{v}$ in the Jacobi--Perron, Brun and in the Selmer algorithm.

In the case of Jacobi--Perron algorithm, the last component of the represented vector has to be the largest component. Therefore, we put $P = \begin{pmatrix}
	0&0&1\\
	0&1&0\\
	1&0&0
\end{pmatrix}$ and find the expansion of $P\vec{v}$. (This means that in this case the matrix of repetend equals to $\pm P(M_1^{m_1} M_2^{m_2})P^{-1}$ for some $m_1,m_2 \in \Z$.) We obtain an eventually periodic expansion $P\vec{v} = \overline{M_{JPA}}$ where
\[
M_{JPA} = A_{JP,1,1}A_{JP,2,4}A_{JP,0,4}A_{JP,0,5}A_{JP,2,4}^{-1}A_{JP,1,1}^{-1} = \begin{pmatrix}
	3 & 9 & 4 \\
	4 & 11 & 5 \\
	5 & 14 & 6
\end{pmatrix} =PM_1M_2^{-3}P^{-1} .
\]

In the Brun algorithm, the vector $\vec{v}$ has a purely periodic expansion equal to $\vec{v} = \overline{M_B}$ where
\[
M_B  = T_{12}T_{23}T_{31}^3T_{12}T_{23}^3T_{31}T_{12}^2 =  \begin{pmatrix}
	20 & 45& 16\\
	16&36&13\\
	13&29&10
\end{pmatrix} = M_1^3 M_2^{-3}.
\]

In the Selmer algorithm, the vector $\vec{v}$ has a purely periodic expansion too. In this case, the expansion is equal to $\vec{v} = \overline{M_S}$ where
\[
M_S = T_{13}T_{21}T_{31}T_{23}T_{12}T_{32} = \begin{pmatrix}
	2&3&1\\
	1&3&1\\
	1&2&1
\end{pmatrix} =  M_2^{-2}.
\]

This means that for the Jacobi--Perron algorithm, we have $m_1 = 1,m_2 = -3$, for the Brun algorithm, we have $m_1 = 3$, $m_2 = -3$ and for the Selmer algorithm we have $m_1 = 0$ and $m_2 = -2$.

\end{example}

The next step is to see if a candidate matrix is in fact a matrix of repetend.
	
\subsection{Decomposition of the candidates on the matrix of repetend} \label{sec:problemb}

	After we have all the candidates on the matrix of repetend of $\vec{v}$, we need to find whether there exists a candidate $M$ on the matrix of repetend for which we can find matrices $R$ and $N$ such that $M = RNR^{-1}$ and such that $\vec{v} = R\overline{N}$.
	In other words, both $R$ and $N$ need to have a decomposition into the matrices from $\A$ and $R\overline{N}$ needs to be formed from an expansion produced by \Cref{alg:MCF_alg} in the given MCF algorithm.
	
	We sum it up in the following proposition.
		
	\begin{proposition}\label{thm:periodic_expansion}
		Let $\vec{v} = \begin{pmatrix}
			y_1\\
			\vdots\\
			y_n
		\end{pmatrix}\in \R^{n}_+$ be a basis of some number field as a vector space over $\Q$.
		
		The vector $\vec{v}$ has an eventually periodic expansion in a unimodular MCF algorithm if and only if there exists $\varepsilon \in U(\Ocal_{\Q( y_1,\dots,  y_n)})$ such that
		$
		{T^{\vec{v}}_\varepsilon}^T 
		$
		(matrix of the linear transformation $t_\varepsilon$ (defined by \eqref{eq:t_alpha}) in the basis $\vec{v}$)
		has a decomposition which is equal to an expansion produced by \Cref{alg:MCF_alg} in the given MCF algorithm. 
	\end{proposition}
	The claim of the proposition is equivalent to the existence of $R,N \in \N^{n,n}$ such that ${T^{\vec{v}}_\varepsilon}^T = RNR^{-1}$ and the matrices $R$ and $N$ have decompositions $R = R_1\cdots R_q$ and $N = N_1\cdots N_p$ such that $\vec{v} = (R_1,\dots,R_q,\overline{N_1,\dots,N_p})$ in the given MCF algorithm.
\begin{proof}
	If the vector $\vec{v}$ has an eventually periodic expansion in a unimodular MCF algorithm, then the rest follows by \Cref{thm:candidates_onMP_algebraic_new} and by the definition of \Cref{alg:MCF_alg}.
	
	To show the converse, let $\varepsilon \in U(\Ocal_{\Q( y_1,\dots,  y_n)})$ be such that $
	{T^{\vec{v}}_\varepsilon}^T
	$ has a decomposition which is equal to an expansion produced by \Cref{alg:MCF_alg}. Let $R,N \in \N^{n,n}$ be such that $
	{T^{\vec{v}}_\varepsilon}^T = RNR^{-1}
	$ and $\vec{v} = RN\dots$ is the expansion of the vector $\vec{v}$ in the given MCF algorithm. Using \Cref{thm:periodicity_of_MCF}, we obtain that $\lambda R^{-1} \vec{v} = N (R^{-1}\vec{v})$ for some $\lambda \in \R$. The matrices $R$ and $N$ are nonnegative and therefore $\lambda >0$.
	Using \Cref{rem:renormalisation}, we obtain that $\vec{v} = R\overline{N}$.
\end{proof}

 All of the well-known algorithms can be defined in a way in which the set $\A$ is equal to the set of transvections or multiples of transvections. 
 Therefore, we can decompose every integer matrix into a product of transvections using \Cref{prop:gen_SL}.
There are numerous decompositions available, and the challenge lies in identifying the specific one that matches the expansion generated by \Cref{alg:MCF_alg} within the given MCF algorithm, if such an expansion exists.
 	In general, this seems to be a difficult question, and it remains to be an open problem for now.

\section{Repetend matrix form and construction of expansions}
\label{sec:application}

As an alternative approach to investigating the possibility to decompose the candidates on the matrix of repetend, we show, on an example, how we can use the knowledge of a matrix $T_\lambda^{\vec{v}}$ for some $\lambda \in \Q(v_1,\dots,v_n)$ of degree $n$, where $\vec{v} = \begin{pmatrix}
	v_1\\
	\vdots\\
	v_n
\end{pmatrix}$, for the construction of expansions of a parametric class of vectors.
We refer to this construction shortly as \emph{repetend matrix form} of the algorithm.

In our example, we construct expansions in the Algebraic Jacobi--Perron algorithm.

The Algebraic Jacobi--Perron algorithm (AJPA) is a Jacobi--Perron type algorithm which was introduced by Tamura and Yasutomi in 2009 (\cite{tamura2009new, tamura2012some}).

For our purposes, we use the AJPA in its homogenous form and study only the dimension $n -1= 2$.

Let $K$ be a cubic number field; $N(v)$ below denotes the norm $N_{K|\Q}(v)$.

Let $\vec{v} = \tvek{v_1}{v_2}{v_0} \in (K\cap\R_+)^{3}$. We have
\[
\mathcal{I}_{AJPA} = \{I_{1,j,k}, I_{2,j,k},I_{3,j,k}\colon j,k \in \N_3\} \text{ with }
\]
\[
I_{1,j,k} = 
\left \{ \tvek{v_1}{v_2}{v_0} \colon \begin{array}{c}\left \lfloor\frac{v_2}{v_1} \right \rfloor = j, \left \lfloor\frac{v_0}{v_1} \right \rfloor = k, \\ v_p >v_1, v_p>v_q ,  \frac{v_1}{\sqrt{|N(v_1)|}}>\frac{v_q}{\sqrt{|N(v_q)|}} \end{array} \right \},
\]
where $p = 0, q= 2$ or $p = 2, q = 0$.

\[
I_{2,j,k} =
 \left \{ \tvek{v_1}{v_2}{v_0} \colon \begin{array}{c} \left \lfloor\frac{v_1}{v_2} \right \rfloor = j, \left \lfloor\frac{v_0}{v_2} \right \rfloor = k, \\ v_p >v_2, v_p>v_q,  \frac{v_2}{\sqrt{|N(v_2)|}}>\frac{v_q}{\sqrt{|N(v_q)|}} \end{array} \right \},
\]
where $p = 0, q= 1$ or $p = 1, q = 0$.

\[
I_{3,j,k} = 
\left \{ \tvek{v_1}{v_2}{v_0} \colon  \begin{array}{c} \left \lfloor\frac{v_1}{v_0} \right \rfloor = j,\left \lfloor\frac{v_2}{v_0} \right \rfloor = k, \\ v_p >v_0, v_p>v_q, \frac{v_0}{\sqrt{|N(v_0)|}}>\frac{v_q}{\sqrt{|N(v_q)|}} \end{array} \right \},
\]
where $p = 2, q= 1$ or $p = 1, q = 2$.

The elements of $\mathcal{I}_{AJPA}$ are pairwise disjoint since all the inequalities in the definition of the sets $I_{i,j,k}$ are strict.

Note that the inequalities defining the intervals above are homogeneous in the sense that their validity does not change when we replace the vector $\tvek{v_1}{v_2}{v_0}$ by its multiple $\tvek{\alpha v_1}{\alpha v_2}{\alpha v_0}$ for any $\alpha\in K\cap\R_+$.

Moreover, let
\[
\A_{AJPA} = \{A_{1,j,k} = T_{21}^jT_{31}^k,A_{2,j,k} = T_{12}^{j}T_{32}^k, A_{3,j,k} = T_{13}^jT_{23}^k \colon j,k \in \N_3\}, 
\]
where the matrices $T_{jk}$ ($j,k \in \{1,\dots,n\}$) are transvections (defined in \Cref{sec:transvections}).

Therefore, the $i$-th step of the algorithm works as follows. If $\vec{v}^{(i)} \in I_{1,j,k}$, then
\[
(v_1^{(i)},v_2^{(i)},v_0^{(i)})^T \mapsto \left (v_1^{(i)}, v_2^{(i)} - j v_1^{(i)}, v_0^{(i)} - k v_1^{(i)}  \right)^T = (v_1^{(i+1)}, v_2^{(i+1)}, v_0^{(i+1)}),
\]
and analogously for $\vec{v}^i$ in other intervals.

\begin{definition}
	The homogenous \emph{AJPA algorithm} is the $(\mathcal{I}_{AJPA},\A_{AJPA})$ MCF algorithm.
	
	The homogenous \emph{AJPA expansion} of a vector is the $(\mathcal{I}_{AJPA},\A_{AJPA})$ MCF expansion.
\end{definition}

Tamura and Yasutomi~\cite[Theorem 2.4]{tamura2009new} give a class of vectors and show they have eventually periodic AJPA expansion.
In the following theorem, we give a larger class and use the results of~\Cref{sec:extracting_information} to prove that the elements of this class have eventually periodic AJPA expansion.

\begin{theorem} \label{thm:AJPA_periodic_class}
	Let $y$ be a positive cubic number such that $-1+ty+sy^2+y^3 = 0$ for some $s,t \in \Q_+$ for which $t>s$ and $ t>\frac{s^2}{4}$. Moreover, let $\vec{v}_{s,t,f,r} = \begin{pmatrix}
		y^2+fy+r\\
		y\\
		1
	\end{pmatrix} $ for some $f,r \in \Z_+ \cup \{0\}$, where $f$ is such that $y^2+fy<y\sqrt{|f^3-sf^2+ft+1|}$ and $ y^2+fy <1$ (especially $f = 0$ fulfils this condition for every $s,t$). The AJPA expansion of $\vec{v}_{s,t,f,r}$ is 
\[
\vec{v}_{s,t,f,r} = A_{3,r,0}A_{2,f,t}\overline{A_{1,t,s}A_{3,t,s}A_{2,s,t}}
\]
for $r>0$ and
\[
\vec{v}_{s,t,f,0} = A_{2,f,t}\overline{A_{1,t,s}A_{3,t,s}A_{2,s,t}}.
\]
\end{theorem}

\begin{proof}
	First of all, we show that $\Q(y)$ is a cubic complex number field. Since $y$ is a cubic number, it is enough to show that the monic minimal polynomial $-1+tx+sx^2+x^3$ of $y$ has a non-real complex root. This happens if the discriminant $\Delta$ of this polynomial is negative. We compute the discriminant:
	\[
	\Delta = -18st+4s^3 + s^2t^2-4t^3-27 = (s^2-4t)(t^2+4s)-2st-27\overset{t>\frac{s^2}{4}}{<}0.
	\]
	Therefore, this equation has one purely real and two complex solutions.
	
	Now we compute the expansion of $\vec{v}_{s,t,f,r}$.
	The coefficients $s,t$ are positive and therefore $y<1$. Moreover, the constant coefficient of the monic minimal polynomial of $y$ is $-1$ and hence $|N(y)| = 1$.
	We start with the case $r>0$. In this case, we have

	\[
	 y^2+fy+r>1>y \quad \text{ and } \quad 1= \frac{1}{|N(1)|}>\frac{y}{|N(y)|} = y \quad \text{ and }  
	\]
	\[
	\lfloor y^2+fy + r\rfloor = r \quad \text{ and } \quad  \lfloor y \rfloor  = 0.
	\]
	Together, we obtain that $\vec{v}_{s,t,f,r} \in I_{3,r,0}$ and $\vec{v}_{s,t,f,r}^{(1)} = \begin{pmatrix}
		y^2+fy\\
		y\\
		1
	\end{pmatrix}$. At this point, we notice, that $\vec{v}_{s,t,f,r}^{(1)} = \vec{v}_{s,t,f,0}$.

Now, we find a matrix $M_0 = {T_\lambda^{\vec{v}_{s,t,f,0}}}^T$ for some cubic number $\lambda \in \Q(y)$. We choose $\lambda$ for which we have $(M_0)_{\bullet,1} = \begin{pmatrix}
	1\\
	1\\
	1
\end{pmatrix}$.
We can do this choice using the explicit connection between this column of the matrix $M_0$ and the coordinates of $\lambda$ in the basis $(y^2+fy,y,1)$ given by \eqref{eq:varepsilon_in_basis} and \eqref{eq:eps1_to_first_column}.

Let $\Id$ be the identity matrix of dimension $3$. Using \Cref{prop:T_epsilon_cubic_case}, we obtain that

\[
Q_{1,\vec{v}_{s,t,f,0}} =
\left(\Id , \begin{pmatrix}
	0 & -f^{2} + f s - t & 1 \\
	1 & -2 \, f + s & 0 \\
	0 & 1 & -f + s
\end{pmatrix}, \begin{pmatrix}
	0 & 1 & f \\
	0 & 0 & 1 \\
	1 & -f + s & t
\end{pmatrix}\right)
\]
and by \Cref{th:weak_exists_Q} we have
\[
M_0 = \begin{pmatrix}
	1 & -f^{2} + f s - t + 1 & f + 1 \\
	1 & -2 \, f + s + 1 & 1 \\
	1 & -f + s + 1 & -f + s + t + 1
\end{pmatrix}.
\]

Now, we use \Cref{lem:Qi_from_M_repetend} and find the triplets $Q_{2,\vec{v}_{s,t,f,0}}$ and $Q_{3,\vec{v}_{s,t,f,0}}$. 
We obtain
{
	\renewcommand{\arraystretch}{1.3}
	\begin{align*}
		& Q_{2,\vec{v}_{s,t,f,0}} = \\
		&
		\left(\begin{pmatrix}
			-\frac{2 \, f^{2} - 3 \, f s + s^{2}}{k} & 1 & -\frac{2 \, f - s}{k} \\
			-\frac{f - s}{k} & 0 & -\frac{1}{k} \\
			-\frac{1}{k} & 0 & -\frac{f^{2} - f s + t}{k}
		\end{pmatrix}, \Id , \begin{pmatrix}
			-\frac{2 \, f - s}{k} & 0 & -\frac{f^{3} - f^{2} s + f t + 1}{k} \\
			-\frac{1}{k} & 0 & -\frac{f^{2} - f s + t}{k} \\
			-\frac{f^{2} - f s + t}{k} & 1 & -\frac{{\left(f^{2} - f s\right)} t + t^{2} + f}{k}
		\end{pmatrix}\right) \text{ and }
		\\
		& Q_{3,\vec{v}_{s,t,f,0}} = \\
		& \left(\begin{pmatrix}
			f - s & -f^{2} + f s - t & 1 \\
			1 & -f & 0 \\
			0 & 1 & 0
		\end{pmatrix}, \begin{pmatrix}
			-f^{2} + f s - t & f^{3} - f^{2} s + f t + 1 & 0 \\
			-f & f^{2} - t & 1 \\
			1 & -2 \, f + s & 0
		\end{pmatrix}, \Id \right),
\end{align*}}
where $k = f^{3} - 2 \, f^{2} s + f s^{2} + {\left(f - s\right)} t - 1$.

We have $1>y>y^2$.
We use \Cref{prop:norms_of_elements} to compute that $|N(y^2+fy)| = f^3-sf^2+ft+1$ and therefore $y = \frac{y}{|N(y)|} >\frac{y^2+fy}{|N(y^2+fy)|}$ by the assumption on $f$. We know that $\lfloor \frac{y^2+fy}{y}\rfloor = f$ so it remains to compute $\left\lfloor \frac{1}{y}\right\rfloor$ in this step of the algorithm. For this reason, we put $M_{0,c,f} = T_{32}^{-c}T_{12}^{-f}M_0T_{12}^{f} T_{32}^c$ and $\vec{v}_{s,t,f,0}^{(0,c)} = T_{32}^{-c}T_{12}^{-f}\vec{v}_{s,t,f,0}$, compute $Q_{1,\vec{v}_{s,t,f,0}^{(0,c)}}$ and use \Cref{prop:comparison_of_elements} to determine the maximum $\widehat{c} \in \Z_+$ such that $\left(\vec{v}_{s,t,f,0}^{(0,c)}\right)_{2}< \left(\vec{v}_{s,t,f,0}^{(0,c)}\right)_{3}$ for all $c \in \Z_+, c\leq \widehat{c}$. We obtain that
\[
\left(\vec{v}_{s,t,f,0}^{(0,c)}\right)_{2}< \left(\vec{v}_{s,t,f,0}^{(0,c)}\right)_{3}
\]
\[ \iff
\]
\[
\left(c^{3} - c^{2} t + 3 \, c^{2} - c s - 2 \, c t + 3 \, c - s - t\right)\left(c^{3} - c^{2} t - c s - 1\right) >0
\]
\[
\iff
\]
\[
((c-t+1)(c^2 +2c) + (c-t)-cs-s)\left(c^{2}(c-t)  - c s - 1\right) >0
\]
and this holds for every $c \leq t-1$. For $c = t$, we obtain
\[
\left(\vec{v}_{s,t,0}^{(0,t)}\right)_{2}< \left(\vec{v}_{s,t,0}^{(0,t)}\right)_{3}\iff (t^2+2t-ts-s)(-ts-1)>0
\]
which does not hold since $t>s$. Therefore $\left \lfloor \frac{1}{y}\right\rfloor = t$.

This means that $\vec{v}_{s,t,f,0} \in I_{2,f,t}$, $\vec{v}_{s,t,f,0}^{(1)} = T_{32}^{-t}T_{12}^{-f}\vec{v}_{s,t,f,0}$ and we put $M_1 = T_{32}^{-t}T_{12}^{-f}M_0T_{12}^{f}T_{32}^t$. From the knowledge of the matrix $M_1$ and by \Cref{lem:Qi_from_M_repetend} we obtain the triplets $Q_{1,\vec{v}_{s,t,f,0}^{(1)}},Q_{2,\vec{v}_{s,t,f,0}^{(1)}}$ and $Q_{3,\vec{v}_{s,t,f,0}^{(1)}}$:

\[
Q_{1,\vec{v}_{s,t,f,0}^{(1)} }= \left(\left(\begin{array}{rrr}
	1 & 0 & 0 \\
	0 & 1 & 0 \\
	0 & 0 & 1
\end{array}\right), \left(\begin{array}{rrr}
	0 & t & 1 \\
	1 & t^{2} + s & t \\
	0 & s t + 1 & s
\end{array}\right), \left(\begin{array}{rrr}
	0 & 1 & 0 \\
	0 & t & 1 \\
	1 & s & 0
\end{array}\right)\right),
\]
\[
Q_{2,\vec{v}_{s,t,f,0}^{(1)} }= \left(\left(\begin{array}{rrr}
	s^{2} - t & 1 & -s \\
	-s & 0 & 1 \\
	s t + 1 & 0 & -t
\end{array}\right), \left(\begin{array}{rrr}
	1 & 0 & 0 \\
	0 & 1 & 0 \\
	0 & 0 & 1
\end{array}\right), \left(\begin{array}{rrr}
	-s & 0 & 1 \\
	1 & 0 & 0 \\
	-t & 1 & 0
\end{array}\right)\right),
\]
\[
Q_{3,\vec{v}_{s,t,f,0}^{(1)} }= \left(\left(\begin{array}{rrr}
	-s & 0 & 1 \\
	1 & 0 & 0 \\
	-t & 1 & 0
\end{array}\right), \left(\begin{array}{rrr}
	0 & 1 & 0 \\
	0 & t & 1 \\
	1 & s & 0
\end{array}\right), \left(\begin{array}{rrr}
	1 & 0 & 0 \\
	0 & 1 & 0 \\
	0 & 0 & 1
\end{array}\right)\right).
\]

Using these triplets and the \Cref{prop:norms_of_elements}, we obtain that $\sqrt{\frac{\left|N\left(\left(\vec{v}_{s,t,f,0}^{(1)}\right)_{3}\right)\right|}{\left|N\left(\left(\vec{v}_{s,t,f,0}^{(1)}\right)_{1}\right)\right|}} = {\sqrt{st+1}}$.

Now, we put $\vec{v}_{s,t,f,0}^{(1,0,b)} = T_{31}^{-b}\vec{v}_{s,t,f,0}^{(1)}$ and $\vec{v}_{s,t,f,0}^{(1,c,0)} = T_{21}^{-c}\vec{v}_{s,t,f,0}^{(1)}$ for $b,c \in \Z_+$.

Again, we put $M_{(1,0,b)} = T_{31}^{-b}M_1T_{31}^b$, we find $Q_{2,\vec{v}_{s,t,f,0}^{(1,0,b)}}$ and use \Cref{prop:comparison_of_elements}. We obtain that

\[
\left(\vec{v}_{s,t,f,0}^{(1,0,b)}\right)_{1}< \left(\vec{v}_{s,t,f,0}^{(1,0,b)}\right)_{3}
\]
\[
\iff
\]
\[
b^{3} + {\left(b + 1\right)} s^{2} + 3 \, b^{2} - 2 \, {\left(b^{2} + 2 \, b + 1\right)} s + {\left(b - s + 1\right)} t + 3 \, b<0
\]
\[
\iff
\]
\[
(b-s+1)(b^2+2b+1-s(b+1)+t)-1<0
\]
which holds if (but not only if) $b\leq s-1$ and $t>(s-b-1)(b+1)$. The second inequality holds if $ b\leq s-1$ and $t >\frac{s^2}{4}$. On the other hand, for $b =s$ we obtain that
\[
\left(\vec{v}_{s,t,f,0}^{(1,0,s)}\right)_{1}< \left(\vec{v}_{s,t,f,0}^{(1,0,s)}\right)_{3} \iff s^2+2s+1-s^2-s+t-1<0 \iff s+t<0
\]
which does not hold.
 It follows that $\left \lfloor \frac{\left(\vec{v}_{s,t,f,0}^{(1)}\right)_{3}}{\left(\vec{v}_{s,t,f,0}^{(1)}\right)_{1}}\right \rfloor = s$.
Similarly, we obtain that
\[
\left(\vec{v}_{s,t,f,0}^{(1,c,0)}\right)_{1}< \left(\vec{v}_{s,t,f,0}^{(1,c,0)}\right)_{2}
\]
\[
\iff
\]
\[
(-c^3 + c^2 t + c s + 1) (-c^3 + c^2 t - 3 c^2 + c s + 2 c t - 3 c + s + t)>0
\]
\[
\iff
\]
\[
((t-c)c^2 + c s + 1) ((c+1)^2(t-c-1)+cs+s+1)>0
\]
which holds for all $c \leq t-1$ but does not hold for $c = t$. This means that $\left \lfloor\frac{\left(\vec{v}_{s,t,f,0}^{(1)}\right)_{2}}{\left(\vec{v}_{s,t,f,0}^{(1)}\right)_{1}}\right\rfloor = t$. Therefore $\left(\vec{v}_{s,t,f,0}^{(1)}\right)_{2}>\left(\vec{v}_{s,t,f,0}^{(1)}\right)_{3}>\left(\vec{v}_{s,t,f,0}^{(1)}\right)_{1}$, $\frac{\left(\vec{v}_{s,t,f,0}^{(1)}\right)_{3}}{\left(\vec{v}_{s,t,f,0}^{(1)}\right)_{1}} = s < \sqrt{st+1} = \sqrt{\frac{\left|N\left(\left(\vec{v}_{s,t,f,0}^{(1)}\right)_{3}\right)\right|}{\left|N\left(\left(\vec{v}_{s,t,f,0}^{(1)}\right)_{1}\right)\right|}}$. This means that $\vec{v}_{s,t,f,0}^{(1)}\in I_{1,t,s}$.

Now, we put $\vec{v}_{s,t,f,0}^{(2)} = T_{21}^{-t}T_{31}^{-s}\vec{v}_{s,t,f,0}^{(1)}$ and $M_2 = T_{21}^{-t}T_{31}^{-s} M_1 T_{31}^{s}T_{21}^{t}$

Again, we use \Cref{lem:Qi_from_M_repetend} to find the triplets $Q_{1,\vec{v}_{s,t,f,0}^{(2)}},Q_{2,\vec{v}_{s,t,f,0}^{(2)}}$ and $Q_{3,\vec{v}_{s,t,f,0}^{(2)}}$. We obtain

\[
Q_{1,\vec{v}_{s,t,f,0}^{(2)} }= \left(\left(\begin{array}{rrr}
	1 & 0 & 0 \\
	0 & 1 & 0 \\
	0 & 0 & 1
\end{array}\right), \left(\begin{array}{rrr}
	0 & 0 & 1 \\
	1 & 0 & -t \\
	0 & 1 & -s
\end{array}\right), \left(\begin{array}{rrr}
	0 & 1 & -s \\
	0 & -t & s t + 1 \\
	1 & -s & s^{2} - t
\end{array}\right)\right),
\]
\[
Q_{2,\vec{v}_{s,t,f,0}^{(2)} }= \left(\left(\begin{array}{rrr}
	t & 1 & 0 \\
	s & 0 & 1 \\
	1 & 0 & 0
\end{array}\right), \left(\begin{array}{rrr}
	1 & 0 & 0 \\
	0 & 1 & 0 \\
	0 & 0 & 1
\end{array}\right), \left(\begin{array}{rrr}
	0 & 0 & 1 \\
	1 & 0 & -t \\
	0 & 1 & -s
\end{array}\right)\right),
\]
\[
Q_{3,\vec{v}_{s,t,f,0}^{(2)} }= \left(\left(\begin{array}{rrr}
	t^{2} + s & t & 1 \\
	s t + 1 & s & 0 \\
	t & 1 & 0
\end{array}\right), \left(\begin{array}{rrr}
	t & 1 & 0 \\
	s & 0 & 1 \\
	1 & 0 & 0
\end{array}\right), \left(\begin{array}{rrr}
	1 & 0 & 0 \\
	0 & 1 & 0 \\
	0 & 0 & 1
\end{array}\right)\right).
\]
Now, we notice that the connection between the triplets $Q_{1,\vec{v}_{s,t,f,0}^{(2)}}$ and $Q_{2,\vec{v}_{s,t,f,0}^{(1)}}$ is the same as in \eqref{eq:per_Q} with matrix of permutation $\widetilde{P}= \begin{pmatrix}
	0 & 1 & 0 \\ 0 & 0 & 1 \\ 1 & 0 & 0
\end{pmatrix}$. Moreover, an analogous connection (with the same matrix of permutation) holds also  between the two pairs of triplets $Q_{2,\vec{v}_{s,t,f,0}^{(2)}}$ and $Q_{3,\vec{v}_{s,t,f,0}^{(1)}}$, $Q_{3,\vec{v}_{s,t,f,0}^{(2)}}$ and $Q_{1,\vec{v}_{s,t,f,0}^{(1)}}$.
Using this fact (analogously as in \Cref{prop:permutations}), we obtain that $\widetilde{P} \vec{v}_{s,t,f,0}^{(1)} = \lambda \vec{v}_{s,t,f,0}^{(2)}$ for some $\lambda >0$ and therefore
\[
\vec{v}_{s,t,f,0} = A_{2,f,t}\overline{A_{1,t,s}A_{3,t,s}A_{2,s,t}}
\]
and
\[
\vec{v}_{s,t,f,r} = A_{3,r,0}A_{2,f,t}\overline{A_{1,t,s}A_{3,t,s}A_{2,s,t}}
\]
for $r>0$.
This proves the claim.
\end{proof}

We now state Theorem 2.4 of \cite{tamura2009new} and give a proof using the last theorem to demonstrate that it indeed covers the class studied in~\cite{tamura2009new}.

Notice that the expansions (and the lengths of periods) of $\vec{v}(m)$ above and the expansions given in Theorem 2.4 in \cite{tamura2009new} slightly differ.
This is due to the fact that we use the homogenous form of the AJPA whereas Tamura and Yatusomi use the non-homogenous form of the AJPA.
The two forms are equivalent and one may transform the expansions from one form to another.

\begin{corollary}[Theorem 2.4 in \cite{tamura2009new}]
		Let $m$ be a positive integer and $\vec{v}(m) = \begin{pmatrix}
			\sqrt[3]{(m^3+1)^2}-m^2\\
			\sqrt[3]{m^3+1}-m\\
			1
		\end{pmatrix} $. All the vectors $\vec{v}(m)$ have eventually periodic homogenous AJPA expansion. The length of the period of $\vec{v}(m)$ is $3$ for every $m >1$ and the length of the period of $\vec{v}(1)$ is $6$.
\end{corollary}
\begin{proof}
	First, we compute explicitly the expansion of $\vec{v}(1)$ and verify its periodicity. We obtain
		\[
	\vec{v}(1) = \overline{A_{1,0,1} A_{2,2,1}A_{3,0,1}A_{1,1,2}A_{2,1,0}A_{3,1,2}}.
	\]
	 Now suppose that $m\geq 2$.
We have $\vec{v}(m) = \begin{pmatrix}
	y^2+2my+m^2\\
	y\\
	1
\end{pmatrix}$ where $y$ is the only real root of the polynomial $g(x) = x^3+3mx^2+3m^2x-1$. We put $s = 3m, t = 3m^2, f = 2m,r = m^2$. 
The function $g(x)$ is increasing and $g(y) = 0$. Therefore, $y^2+2my<1$ and $y^2<y$. Now, we show that $y^2 + fy < \sqrt{|f^3-sf^2+ft+1|}y$. This is equivalent to $ y^2 + 2my <\sqrt{|2m^3+1|}y$. For $m\geq 3$ it follows from the fact that $y^2<y$.
For $m = 2$, this condition is equivalent to $y<\sqrt{17}-4$. We compute that $g(\sqrt{17}-4)>0$ which implies (by the monotony of $g$) that $y<\sqrt{17}-4$. This means that $\vec{v}(m)$ fulfils the assumptions of \Cref{thm:AJPA_periodic_class} for all $m\geq 2$ and therefore
\[
\vec{v}(m) = A_{3,m^2,0}A_{2,2m,3m^2}\overline{A_{1,3m^2,3m}A_{3,3m^2,3m}A_{2,3m,3m^2}}
\]
for all $m\geq 2$.
\end{proof}

\section{Conclusion}

Let us conclude with several remarks and further research directions.

The case of vectors from a totally real number field is of special interest when considering MCF expansions.
While our~\Cref{prop:purly_per} showed that a vector with all coordinates being totally positive cannot have a purely periodic expansion, quite a few examples of eventually periodic expansions are known \cite{KST}.

In particular, \cite{KST} explicitly considered the connection of Jacobi--Perron expansions with universal quadratic forms (and suitable small, ``indecomposable'' elements in the number field) already mentioned in the Introduction. Although their results are promising, they remain only partial and suggest that the more general approach outlined in our present paper may be needed in order to obtain a tight connection between MCFs and indecomposables. Specifically, can one find a suitable decomposition of the candidate for the matrix of the repetend that would yield indecomposables in the form of certain ``(semi-)convergents'' to the expansion?

As we demonstrated in~\Cref{sec:application}, the repetend matrix form of algorithms can be useful computationally, as it avoids working with small real numbers avoiding potential precision issues and at the same time, we can easily find the norms of the components of the represented vector. This should be very convenient in practical computer implementations.
Note that the repetend matrix form is not limited to $n=3$ and Algebraic Jacobi--Perron algorithm; however, its formal generalization requires further study.

Finally, the holy grail in the area of MCFs is establishing that some vectors do not have eventually periodic expansions. For example, computational evidence \cite{Voutier} suggests that this is the case for the Jacobi--Perron expansion of $(1,\sqrt[3]{4},\sqrt[3]{4^2})$. In fact, Voutier [personal communication] conjectured that the positive integers $m$ with eventually periodic JPA expansion of $(1,\sqrt[3]{m},\sqrt[3]{m^2})$ have density 0.
These problems are notoriously hard, but the approach outlined in our paper could present a starting point.

\section{Acknowledgements}

We thank the anonymous referee for several very helpful suggestions.

H. Řada was supported by the Grant Agency of the Czech Technical University in Prague,
grant No. SGS20/183/OHK4/3T/14.
V. Kala was supported by Czech Science Foundation (GA\v CR) grant 21-00420M and by Charles University Research Centre program UNCE/SCI/022.
Š. Starosta acknowledges support of the OP VVV MEYS funded project CZ.02.1.01/0.0/0.0/16\_019/0000765 ``Research Center for Informatics''.

The computer experiments were done using the computer algebra system SageMath \cite{sage_2022}.

	\bibliographystyle{siam}
	\IfFileExists{biblio.bib}{\bibliography{biblio}}{\bibliography{biblio/biblio}}

\begin{thebibliography}{10}

\bibitem{Hermite1850_2}
{\em Extraits de lettres de {M}. {Ch}. {H}ermite à {M}. {J}acobi sur
  différents objets de la théorie des nombres}, J. Reine Angew. Math., 40
  (1850), pp.~279--315.
\newblock Deuxi\`eme lettre.

\bibitem{abrate2013periodic}
{\sc M.~Abrate, S.~Barbero, U.~Cerruti, and N.~Murru}, {\em Periodic
  representations for cubic irrationalities}, Fibonacci Quart., 50 (2012),
  pp.~252--264.

\bibitem{alaca2004introductory}
{\sc {\c{S}}.~Alaca and K.~S. Williams}, {\em Introductory algebraic number
  theory}, Cambridge University Press, 2004.

\bibitem{barbeau2006pell}
{\sc E.~J. Barbeau}, {\em Pell’s equation}, Springer Science \& Business
  Media, 2006.

\bibitem{bernstein1964periodical}
{\sc L.~Bernstein}, {\em Periodical continued fractions for irrationals of
  degree n by {J}acobi's algorithm}, J. Reine Angew. Math., 213 (1964),
  pp.~31--38.

\bibitem{bernstein1964periodicity}
{\sc L.~Bernstein}, {\em Periodicity of {J}acobi's algorithm for a special type
  of cubic irrationals}, J. Reine Angew. Math., 213 (1964), pp.~137--146.

\bibitem{bernstein1964periodische}
{\sc L.~Bernstein}, {\em {Periodische {J}acobische {A}lgorithmen f{\"u}r eine
  unendliche {K}lasse algebraischer {I}rrationalzahlen vom {G}rade n und einige
  unendliche {K}lassen kubischer {I}rrationalzahlen.}}, J. Reine Angew. Math.,
  214 (1964), pp.~76--83.

\bibitem{bernstein1965rational}
\leavevmode\vrule height 2pt depth -1.6pt width 23pt, {\em Rational
  approximations of algebraic irrationals by means of a modified
  {J}acobi-{P}eron algorithm}, Duke Math. J., 32 (1965), pp.~161--176.

\bibitem{blomer2017rank}
{\sc V.~Blomer and V.~Kala}, {\em On the rank of universal quadratic forms over
  real quadratic fields}, Doc. Math., 23 (2018), pp.~15--34.

\bibitem{bouhamza1984algorithme}
{\sc M.~Bouhamza}, {\em Algorithme de {J}acobi-{P}erron dans les corps de
  nombres de degr{\'e} 4}, Acta Arith., 2 (1984), pp.~141--145.

\bibitem{bouhamza1984jacobi}
\leavevmode\vrule height 2pt depth -1.6pt width 23pt, {\em Jacobi {P}erron
  algorithm in cubic number-fields}, B. Sci. Math., 108 (1984), pp.~101--111.

\bibitem{brentjes}
{\sc A.~J. Brentjes}, {\em Multi-dimensional continued fraction algorithms},
  vol.~145 of Mathematical Centre Tracts, Mathematisch Centrum, 1981.

\bibitem{brun1920}
{\sc V.~Brun}, {\em En generalisation av {K}jedebroken}, Skrifter utgit av
  Videnskapsselskapet i Kristiania I Matematisk-naturvidenskabelig klasse, No.
  6 (1920).

\bibitem{Bryuno_1997}
{\sc A.~D. Bryuno and V.~I. Parusnikov}, {\em Comparison of various
  generalizations of continued fractions}, Math. Notes, 61 (1997),
  pp.~278--286.

\bibitem{buchmann1985generalization}
{\sc J.~Buchmann}, {\em A generalization of {V}oronoi's unit algorithm {I}}, J.
  Number Theory, 20 (1985), pp.~177--191.

\bibitem{buchmann1985generalizationII}
\leavevmode\vrule height 2pt depth -1.6pt width 23pt, {\em A generalization of
  {V}oronoi's unit algorithm {II}}, J. Number Theory, 20 (1985), pp.~192--209.

\bibitem{conder1992presentations}
{\sc M.~Conder, E.~Robertson, and P.~Williams}, {\em Presentations for
  3-dimensional special linear groups over integer rings}, Proc. Amer. Math.
  Soc., 115 (1992), pp.~19--26.

\bibitem{cusick1987table}
{\sc T.~W. Cusick and L.~Schoenfeld}, {\em A table of fundamental pairs of
  units in totally real cubic fields}, Math. Comp., 48 (1987), pp.~147--158.

\bibitem{david1957contribution}
{\sc M.~David}, {\em Contribution a l'{\'e}tude algorithmique des
  approximations rationnelles simultan{\'e}es de deux irrationnels.
  {A}pplication au cas cubique. {\'E}tude locale des vari{\'e}t{\'e}s
  analytiques complexes. {I}nvariance de la dimension complexe par les
  transformations pseudo-conformes}, Impr. Durand, 1957.

\bibitem{delone1964theory}
{\sc B.~Delone and D.~Faddeev}, {\em The Theory of Irrationalities of the Third
  Degree}, vol.~10 of Translations of Mathematical Monographs, American
  Mathematical Society, 1964.

\bibitem{dubois1975algorithme}
{\sc E.~Dubois and R.~Paysant-Le~Roux}, {\em Algorithme de {J}acobi-{P}erron
  dans les extensions cubiques}, C. R. Math. Acad. Sci. Paris, 280 (1975),
  pp.~183--186.

\bibitem{fogg2002substitutions}
{\sc N.~P. Fogg}, {\em Substitutions in dynamics, arithmetics and
  combinatorics}, Springer Science \& Business Media, 2002.

\bibitem{greiter1977mehrdimensionale}
{\sc G.~Greiter}, {\em Mehrdimensionale {K}ettenbr{\"u}che}, PhD thesis,
  Technische Universit{\"a}t M{\"u}nchen, 1977.

\bibitem{jacobi}
{\sc E.~Heine and C.~Jacobi}, {\em Allgemeine {T}heorie der
  kettenbruch{\"a}hnlichen {A}lgorithmen, in welchen jede {Z}ahl aus drei
  vorhergehenden gebildet wird.}, J. Reine Angew. Math., 69 (1868), pp.~29--64.

\bibitem{kala2016universal}
{\sc V.~Kala}, {\em Universal quadratic forms and elements of small norm in
  real quadratic fields}, Bull. Aust. Math. Soc., 94 (2016), pp.~7--14.

\bibitem{KST}
{\sc V.~Kala, E.~Sgallová, and M.~Tinková}, {\em Arithmetic of cubic number
  fields: {Jacobi--Perron, Pythagoras,} and indecomposables}, preprint
  available at \url{https://arxiv.org/abs/2303.00485},  (2023).

\bibitem{kala2020universal}
{\sc V.~Kala and M.~Tinková}, {\em {Universal Quadratic Forms, Small Norms,
  and Traces in Families of Number Fields}}, Int. Math. Res. Not., 2023 (2022),
  pp.~7541--7577.

\bibitem{kala2021lifting}
{\sc V.~Kala and P.~Yatsyna}, {\em Lifting problem for universal quadratic
  forms}, Adv. Math., 377 (2021), p.~107497.

\bibitem{karpenkov2021periodic}
{\sc O.~Karpenkov}, {\em On a periodic {J}acobi-{P}erron type algorithm},
  preprint available at \url{https://arxiv.org/abs/2101.12627},  (2021).

\bibitem{karpenkov2013geometry}
\leavevmode\vrule height 2pt depth -1.6pt width 23pt, {\em Geometry of
  continued fractions}, vol.~26 of Algorithms and Computation in Mathematics,
  Springer Berlin, Heidelberg, 2~ed., 2022.

\bibitem{karpenkov2021hermite}
\leavevmode\vrule height 2pt depth -1.6pt width 23pt, {\em On {H}ermite's
  problem, {J}acobi-{P}erron type algorithms, and {D}irichlet groups}, Acta
  Arith. 203,  (2022), pp.~27--48.

\bibitem{levesque1976jacobi}
{\sc C.~Levesque}, {\em {J}acobi-{P}erron algorithms and units in some
  algebraic number fields}, PhD thesis, Illinois Institute of Technology, 1976.

\bibitem{levesque1977class}
\leavevmode\vrule height 2pt depth -1.6pt width 23pt, {\em A class of periodic
  {J}acobi-{P}erron algorithms in pure algebraic number fields of degree n >
  2}, Manuscripta Math., 22 (1977), pp.~235--269.

\bibitem{murru2015periodic}
{\sc N.~Murru}, {\em On the periodic writing of cubic irrationals and a
  generalization of {R}{\'e}dei functions}, Int. J. Number Theory, 11 (2015),
  pp.~779--799.

\bibitem{paysant1975nouvelles}
{\sc R.~Paysant-Le~Roux and E.~Dubois}, {\em De nouvelles classes d'entiers
  pour lesquelles on connaît des développements périodiques par l'algorithme
  de {J}acobi--{P}erron dans $\mathbb{Q}(m^{{1}/{(n+1)}})$}, C. R. Math. Acad.
  Sci. Paris, 280 (1975), pp.~57--59.

\bibitem{Perron}
{\sc O.~Perron}, {\em Grundlagen für eine {T}heorie des {J}acobischen
  {K}ettenalgorithmus}, Math. Ann.,  (1907).

\bibitem{perron1935satz}
{\sc O.~Perron}, {\em Ein {S}atz {\"u}ber {J}acobi-{K}etten zweiter {O}rdnung},
  Ann. Sc. Norm. Super. Pisa Cl. Sci., 4 (1935), pp.~133--138.

\bibitem{poincare}
{\sc H.~Poincar{\'e}}, {\em Sur une g{\'e}n{\'e}ralisation des fractions
  continues}, C. R. Math. Acad. Sci. Paris, 1 (1884), pp.~1014--1016.

\bibitem{raju1976periodic}
{\sc N.~Raju}, {\em Periodic {J}acobi-{P}erron algorithms and fundamental
  units}, Pacific Journal of Mathematics, 64 (1976), pp.~241--251.

\bibitem{schweiger1995}
{\sc F.~Schweiger}, {\em Invariant measures for fully subtractive algorithms},
  Anz. {\"O}sterreich. Akad. Wiss. Math.-Natur. Kl, 131 (1995), pp.~25--30.

\bibitem{schweiger2000}
{\sc F.~Schweiger}, {\em Multidimensional continued fractions}, Oxford
  University Press, 2000.

\bibitem{selmer}
{\sc E.~S. Selmer}, {\em Om flerdimensjonal kjedebr{\o}k}, Nordisk Matematisk
  Tidskrift,  (1961), pp.~37--43.

\bibitem{tamura2009new}
{\sc J.-i. Tamura and S.-i. Yasutomi}, {\em A new multidimensional continued
  fraction algorithm}, Math. Comp., 78 (2009), pp.~2209--2222.

\bibitem{tamura2012some}
{\sc J.-i. Tamura and S.-i. Yasutomi}, {\em Some aspects of multidimensional
  continued fraction algorithms ({F}unctions in number theory and their
  probabilistic aspects)}, RIMS Kôkyûroku Bessatsu, B34 (2012), pp.~463--475.

\bibitem{sage_2022}
{\sc {The Sage Developers}}, {\em {S}ageMath, the {S}age {M}athematics
  {S}oftware {S}ystem ({V}ersion 9.5)}, 2022.
\newblock {\texttt http://www.sagemath.org}.

\bibitem{voronoi1896generalization}
{\sc G.~F. Voronoi}, {\em On a generalization of the algorithm of continued
  fractions}, PhD thesis, Warsaw University, 1896.

\bibitem{Voutier}
{\sc P.~Voutier}, {\em Families of periodic {J}acobi-{P}erron algorithms for
  all period lengths}, J. Number Theory, 168 (2016), pp.~472--486.

\bibitem{yatsyna2019lower}
{\sc P.~Yatsyna}, {\em A lower bound for the rank of a universal quadratic form
  with integer coefficients in a totally real number field}, Comment. Math.
  Helv., 94 (2019), pp.~221--239.

\end{thebibliography}
	
\end{document}